\documentclass[10pt]{amsart}

\usepackage[all]{xypic}

\usepackage{latexsym}

\usepackage{amssymb}
\usepackage{amsfonts}
\usepackage{amscd}
\usepackage{amsmath,amsthm}

\newtheorem{lemma}{Lemma}[section]
\newtheorem{proposition}[lemma]{Proposition}
\newtheorem{theorem}[lemma]{Theorem}
\newtheorem{corollary}[lemma]{Corollary}

{
}

\theoremstyle{definition}

\theoremstyle{remark}

\numberwithin{equation}{section}

\newenvironment{pf}{\noindent{\bf Proof.}}{\hfill $\square$\medskip}

\def\CC{{\mathbb C}}

\def\HH{{\mathbb H}}

\def\MM{{\mathbb M}}
\def\NN{{\mathbb N}}

\def\PP{{\mathbb P}}

\def\ZZ{{\mathbb Z}}

\def\qol{{\bar q}}

\def\alphaol{{\bar{\alpha}}}

\def\phiol{{\bar \phi}}

\def\Mol{{\bar M}}

\def\0ol{{\bar 0}}
\def\1ol{{\bar 1}}
\def\2ol{{\bar 2}}
\def\ol2{{\bar 2}}
\def\3ol{{\bar 3}}
\def\4ol{{\bar 4}}
\def\5ol{{\bar 5}}
\def\6ol{{\bar 6}}
\def\7ol{{\bar 7}}
\def\8ol{{\bar 8}}
\def\9ol{{\bar 9}}

\def\bold0{{\bf 0}}
\def\bold1{{\bf 1}}
\def\bold2{{\bf 2}} 
\def\bold3{{\bf  3}}
\def\bold4{{\bf 4}}
\def\bold5{{\bf 5}}
\def\bold6{{\bf 6}}
\def\bold7{{\bf 7}}
\def\bold8{{\bf 8}}
\def\bold9{{\bf 9}}

\def\P2Skly{\PP^2_{Skly}}

\def\coker{\operatorname {coker}}

\def\End{\operatorname {End}}
\def\Ext{\operatorname {Ext}}

\def\gr{\operatorname {gr}}

\def\Hom{\operatorname {Hom}}

\def\im{\operatorname {im}}

\def\ker{\operatorname {ker}}

\def\MCM{\operatorname {MCM}}

\def\uMCM{\operatorname {\underline{MCM}}}

\def\rad{\operatorname {rad}}

\def\th{\operatorname {th}}    % for writing n^{th}
\def\Tor{\operatorname {Tor}}

\def\depth{\operatorname{depth}}
\def\sD^b{\operatorname{ sD^b}}

\def\dim{\operatorname{dim}}

\def\End{\operatorname{End}}

\def\Ext{\operatorname{Ext}}
\def\Fdim{{\sf Fdim}}
\def\fdim{{\sf fdim}}

\def\Fract{\operatorname{Fract}}

\def\GKdim{\operatorname{GKdim}}
\def\gldim{\operatorname{gldim}}

\def\mod {{\sf mod }}

\def\Gr{{\sf Gr}}
\def\Gr{{\sf Gr}}

\def\gr{{\sf gr}}
\def\H{\operatorname{H}}

\def\Hom{\operatorname{Hom}}

\def\id{\operatorname{id}}

\def\Im{\operatorname{Im}}

\def\liminj{\varinjlim}
\def\Lin{\operatorname{\sf Lin}}

\def\mod{{\sf mod}}
\def\Mod{{\sf Mod}}

\def\pdim{\operatorname{pdim}}

\def\Proj{\operatorname{Proj}}

\def\QGr{\operatorname{\sf QGr}}
\def\qgr{\operatorname{\sf qgr}}

\def\rank{\operatorname{rank}}

\def\RuHom{\operatorname {R\!\uHom}}

\def\Spec{\operatorname{Spec}}

\def\uExt{\operatorname{\underline{Ext}}}

\def\uHom{\operatorname{\underline{Hom}}}

\def\ul1{\operatorname{\underline{1}}}

\def\a{\alpha}
\def\b{\beta}

\def\fm{{\mathfrak m}}

\def\fsl{{\mathfrak s}{\mathfrak l}}

\def\sA{{\sf A}}

\def\sD{{\sf D}}

\def\cal{\mathcal}

\def\cD{{\cal D}}

\def\cF{{\cal F}}
\def\cG{{\cal G}}

\def\cL{{\cal L}}
\def\cM{{\cal M}}

\def\cO{{\cal O}}

\def\coh{{\sf coh}}

\def\Qcoh{{\sf Qcoh}}

\def\dirlim{\mathop{\vtop{\baselineskip -100pt\lineskip -1pt\lineskiplimit 0pt
\setbox0\hbox{lim}\copy0\hbox to \wd0{\rightarrowfill}}}\limits}
\def\invlim{\mathop{\vtop{\baselineskip -100pt\lineskip -1pt\lineskiplimit 0pt
\setbox0\hbox{lim}\copy0\hbox to \wd0{\leftarrowfill}}}\limits}

\def\I11{{1 \kern -0.8pt \! \mbox{l}}} 
\def\mumu{{\mu\kern-4.2pt\mu}}
\def\bfmu{{\mu\kern-4.2pt\mu}}
\def\2slash{\backslash \! \backslash}

\def\boxtimes{\setbox0\hbox{$\Box$}\copy0\kern-\wd0\hbox{$\times$}}

\begin{document}
\pagenumbering{arabic}

\title{Non-commutative quadric surfaces}
\author{S. Paul Smith and M. Van den Bergh}
\address{S. Paul Smith: Department of Mathematics, Box 354350, Univ.  Washington, Seattle, WA 98195}
\address{M. Van den Bergh: Departement WNI, Limburgs Universitair Centrum, Universitaire Campus, Building D, 3590 Diepenbeek, Belgium} 
\email{\tt smith@math.washington.edu,  \hskip .02in michel.vandenbergh@luc.ac.be}



\thanks{S. Paul Smith was supported by the National Science
Foundation under  Grants DMS-0070560 and DMS-024572. M. Van den Bergh is a director of research at the FWO. }

%\tableofcontents

\begin{abstract}
The 4-dimensional Sklyanin algebra is the homogeneous coordinate ring of a noncommutative analogue of projective 3-space. The degree-two component of the algebra contains a 2-dimensional subspace of 
central elements. The zero loci of those central elements, except 0, form a pencil of non-commutative quadric 
surfaces, We show that the behavior of this pencil is similar to that of a generic pencil of quadrics in the commutative projective 3-space. There are exactly four singular quadrics in the pencil. The singular and non-singular quadrics are characterized by whether they have one or two rulings by non-commutative lines. The Picard groups of the smooth quadrics are free abelian of rank two. The alternating sum of dimensions of Ext groups allows us to define an intersection pairing on the Picard group of the smooth noncommutative quadrics. A surprise is that a smooth noncommutative quadric can sometimes contain a ``curve'' having self-intersection number -2.  Many of the methods used in our paper are noncommutative versions of  methods developed by Buchweitz, Eisenbud and Herzog: in particular, the correspondence between the geometry of a
quadric hypersurface and maximal Cohen-Macaulay modules over its homogeneous coordinate ring plays 
a key role. An important aspect of our work is to introduce definitions of non-commutative analogues of the familiar commutative terms used in this abstract. We expect the ideas we develop here for 2-dimensional non-commutative quadric hypersurfaces will apply to higher dimensional non-commutative 
quadric hypersurfaces and we develop them in sufficient generality to make such applications possible. 
\end{abstract}

\maketitle

%-------------------------------------------------------------------------------------------------------------------------------------------

\section{Introduction}

\subsection{}
Many non-commutative analogues of $\PP^3$ contain non-commutative 
quadric hypersurfaces. This paper studies these non-commutative
quadrics and the consequences of their
existence for the ambient non-commutative $\PP^3$. 

For example, we establish a simple ``geometric'' criterion for recognizing 
when a non-commutative quadric surface is smooth: it is smooth if and only if
it has two rulings. Of course, a key point is to define the terms.

The smoothness result allows us to make further comparisons 
with the commutative case. For example, a generic pencil of quadrics 
in $\PP^3$ has exactly four singular members and we show the same 
is true for the pencil of non-commutative quadrics in the 
Sklyanin analogue of $\PP^3$.

\subsection{}
A non-commutative quadric surface $Q$ is defined implicitly by defining a 
Grothendieck category
$\Qcoh \, Q$ that plays the role of quasi-coherent sheaves on it.
We say that $Q$ is smooth of dimension 
two if $\Ext^3_Q(-,-)$ vanishes everywhere on $\Qcoh \, Q$ but  $\Ext^2_Q(-,-)$ does not.

Deciding whether a commutative variety is smooth is a local problem: one examines the local
rings at its points. One can also use the Jacobian criterion on affine patches. Deciding
 whether a non-commutative variety is smooth is a different  kind of problem because 
the variety can have few closed points, sometimes none at all. One cannot 
check smoothness by checking the homological properties of 
 individual points. In this sense, smoothness is not a local
property and global methods must be used. 
In particular, there is no analogue of the 
Jacobian criterion and singular non-commutative quadrics 
need not have a singular point.

Theorem 5.6 shows that the smoothness of a quadric hypersurface in a non-commutative $\PP^3$ (where these 
terms have to be defined appropriately) is equivalent to the semisimplicity of a certain finite dimensional algebra. 

\subsection{}
We will define non-commutative quadric surfaces as degree two
hypersurfaces in suitable non-commutative analogues of $\PP^3$, the
latter being a non-commutative space of the form $\Proj S$ 
where $S$ is a not-necessarily-commutative connected graded ring
having  properties like those of the commutative polynomial
ring in four variables (see section \ref{sect.qu.Pn} for a precise definition). 
Thus $Q=\Proj A$ where $A=S/(z)$ and $z \in S_2$ is a central 
regular element. (If $z$ were a normal regular element we could
replace $S$ by a suitable Zhang twist in which $z$ becomes central,
so there is no loss of generality in assuming $z$ is central.)
Amongst other things, $S$ is required to be a Koszul algebra and
this implies that $S/(z)$ is also Koszul.

\subsection{}
Let $S$ denote a 4-dimensional Sklyanin algebra \cite{Skly}, 
\cite{Sm1}, \cite{SS}, \cite{TV}, \cite{vdB2}. In this case we write
$$
\PP^3_{\sf Skly} = \Proj S.
$$
The common zero locus of the two linearly independent degree-two central elements in $S$ is  
commutative elliptic curve $E$. 
The zero loci of linear combinations of these two central elements form
a pencil of non-commutative quadrics $Q =\Proj S/(z) \subset \PP^3_{\sf Skly}$.
Exactly four of these non-commutative quadrics are singular
(Theorem \ref{thm.sing.skly.qus}). The base locus of the pencil is $E$.
This is a direct analogue of the 
commutative case: the base locus of a generic pencil of quadrics 
in $\PP^3$  is a quartic elliptic curve, and exactly
four members of that pencil are singular. 

\subsection{}
The method we use to understand these non-commutative quadrics follows 
that of Buchweitz, Eisenbud, and Herzog in their paper \cite{BEH} 
on maximal Cohen-Macaulay modules over quadrics, and especially 
the approach in the appendix of their paper. 
As is well known, a quadric hypersurface is smooth if and only if
the even Clifford algebra determined by its defining equation is
semisimple.
The results in \cite{BEH} establish a duality between the maximal Cohen-Macaulay modules over the coordinate ring of the quadric and the derived category of the Clifford algebra. 
We associate to our non-commutative quadrics $Q=\Proj A$ 
finite dimensional algebras $C$ that are analogues of
even Clifford algebras and establish the ``same'' duality. 

\subsection{}
Let $A^!$ be the quadratic dual of $A$.
Because $A$ is a hypersurface ring, $\Proj A^!$ is an affine space.
The algebra $C$ is a coordinate ring of this space in the sense 
that $\Qcoh(\Proj A^!) \cong \Mod C$.
We show that $Q$ is smooth if and only if $C$ is semisimple if and only 
if there are two distinct non-commutative ``rulings'' on $Q$.
We show that the ``lines''
on $Q$ determine  
two-dimensional simple $C$-modules; because the dimension
of $C$ is 8, it is semisimple if it has two non-isomorphic 
two-dimensional simple modules. The method by which we 
associate a $C$-module to a line on $Q$ uses the fact that $A$
is a Koszul algebra, and that the lines on
$Q$ determine graded maximal Cohen-Macaulay $A$-modules.

\subsection{}
Although quantum $\PP^2$s have been classified  and are 
well-understood in some regards, the same is not true for quantum
$\PP^3$s.  The results in this paper are a step towards gaining a 
similar understanding of another class of non-commutative
surfaces. In
sections \ref{sect.lines} and \ref{sect.pts} we obtain good
information about the points and lines on such surfaces.

Furthermore, if the non-commutative quadric $Q$ is smooth there is an isomorphism
$K_0(Q) \cong K_0(\PP^1 \times \PP^1)$ of
Grothendieck groups that is compatible with the Euler forms 
$(-,-)=\sum (-1)^i \dim \Ext^i_Q(-,-)$. More interestingly, 
the effective cones for $Q$ and $\PP^1 \times \PP^1$ need not match 
up under this isomorphism:
sometimes $Q$ contains, in effect, a $-2$-curve.

All unexplained terminology for non-commutative spaces can be found 
in either \cite{SmSub} or \cite{vdB3}.

\section{Preliminaries}

{\it Throughout $k$ denotes a field and $A$ denotes a two-sided noetherian connected graded 
$k$-algebra. }

\medskip

The Hilbert series of a graded $k$-vector $V$ having 
finite-dimensional components is the formal series
$$
H_V(t):=\sum_n (\dim_k V_n) t^n.
$$

\subsection{Graded modules}

The category of graded right $A$-modules with degree zero module homomorphisms is denoted by 
 $\Gr A$  and $\gr  A$ is the full subcategory of $\Gr A$ consisting of noetherian modules. 
We write $\sD^b(\gr A)$, or just $\sD^b(A)$, for the associated bounded derived category.

We write $\Ext^i_{\Gr A}(M,N)$ for the extension groups in 
$\Gr A$, and define
$$
\uExt^*_A(M,N):=\bigoplus_{i \in \ZZ} \Ext^i_{\Gr A}(M,N(i)).
$$

\subsection{Syzygies}

When $M \in \gr A$ we write $\Omega^iM$ for the $i^{\th}$ syzygy in $\gr  A$ 
obtained from a minimal graded resolution 
of $M$.  Since $A$ is connected graded, $\uExt_A^{i}(\Omega^dM,k) \cong \uExt_A^{i+d}(M,k)$ for all $i \ge
0$. We often write $\Omega M$ for $\Omega^1 M$.

\subsection{Linear resolutions}

An $M$ in $\gr A$ has a {\sf linear resolution} if for all $i$ the $i^{\th}$ term in its minimal 
projective resolution is a direct sum of copies of $A(-i)$ or,
equivalently, if  $\uExt_A^i(M,k)_j =0$ whenever $i+j \ne 0$. 
We write $\Lin(A)$ for the full subcategory  of $\gr A$
consisting of modules having a linear resolution;
$\Lin (A)$ is closed under direct summands and extensions. 

If $M \in \Lin (A)$, then $(\Omega^nM)(n) \in \Lin (A)$ too.

If $M \in \Lin (A)$, then
\begin{equation}
\label{fnl.eq}
H_{\uExt^*_A(M,k)}(t)H_A(-t) = H_M(-t).
\end{equation}

\subsection{Koszul duality} 

See \cite{BGS} for basic information about Koszul algebras. 

Let $A$ be a connected Koszul algebra and $A^!$ its quadratic dual. 

The Koszul property says that the natural homomorphism 
$A^! \to \uExt_A^*(k,k)$ is an isomorphism of graded $k$-algebras. 
If $M$ is a graded $A$-module, the Yoneda product makes
$\uExt^*_A(M,k)$ a graded left $A^!$-module with degree $i$
component $\uExt^i_A(M,k)$.

A graded $A$-module $M$ is {\sf stably linear} if $M_{\ge n}(n)$ has a linear
resolution for $n \gg 0$.
We write $\sD^b_{sl}(A)$ for the full subcategory of $\sD^b(A)$ of
complexes having stably linear homology. 

By \cite{SvdB} there is a duality
$$
K:\sD^b_{sl}(A) \to \sD^b_{sl}(A^!)
$$
given by
$$
KM=T\bigl(\RuHom_A(M,k))
$$
where $T$ is the re-grading functor 
$$
(TV)^i_j=V^{i+j}_{-j}
$$
where the upper index is the homological 
degree and the lower index  the grading degree.
The duality $K$ restricts to a duality 
$$
\Lin(A) \to \Lin(A^!), \qquad M \mapsto \uExt^*_A(M,k)
$$
with degree $i$ component $\uExt^i_A(M,k)$.
The Koszul duality functor $K$ satisfies
\begin{equation}
\label{eq.K.twist}
K(M[1])  \cong (KM)[-1]
\qquad
\hbox{and}
\qquad
K(M(1))  \cong (KM)[-1](1).
\end{equation}

\begin{theorem}
[J\"orgensen]
\cite[Thm. 3.1]{J}
Let $A$ be a two-sided noetherian, connected, graded $k$-algebra that is 
Koszul and has a balanced dualizing complex \cite{vdB7}, \cite{YZ}.
Then every finitely generated $A$-module is stably linear.
Thus $$\sD^b_{sl}(A)=\sD^b(A).$$
\end{theorem}

\subsection{Cohen-Macaulay rings and modules}

Let $A$ be a right and left noetherian, connected, graded $k$-algebra
having a balanced dualizing complex $R^\cdot$ \cite{vdB7}, \cite{YZ}.

We say that $A$ is {\sf Cohen-Macaulay} of {\sf depth $d$} if 
there is an $A$-$A$-bimodule $\omega_A$ such that 
$R^\cdot \cong \omega_A[d]$.
We call $\omega_A$ the {\sf dualizing module} for $A$. 
By \cite[Propostion 7.9]{AZ}, $\omega_A$ is finitely generated on
each side.
We say $A$ is {\sf Gorenstein} if it is Cohen-Macaulay and 
$\omega_A$ is an invertible bimodule. This is equivalent to
the requirement that $\omega_A$ is isomorphic to $A(\ell)$ for 
some $\ell$ as both a right and as a left module.

The local cohomology functors 
$$
H^i_{\fm}(-)= \liminj\uExt^i_A(A/A_{\ge n},-)
$$
are defined on graded right $A$-modules.
Here $\fm$ denotes the maximal ideal $A_{\ge 1}$. We write
$H^i_{\fm^\circ}$ for the local cohomology modules for left
modules.
The {\sf depth} of an $A$-module $M$ is the smallest
integer $i$ such that $H^i_{\fm}(M) \ne 0$. 
A finitely generated module $M$ is {\sf Cohen-Macaulay} if 
either $M=0$ or only one $H^i_{\fm}(M)$ is non-zero. 
For the rest of this section we assume that $A$ is Cohen-Macaulay 
of depth $d$ in the sense 
of the previous paragraph. Then $A$ is a Cohen-Macaulay $A$-module 
of depth $d$ in the sense of the present paragraph.
Furthermore, there is an isomorphism
$$
\omega_A \cong H^d_{\fm}(A)^*
$$
of $A$-$A$-bimodules and, for every $M \in \mod  A$,
\begin{equation}
\label{dual}
\uExt^i_A(M,\omega_A) \cong H_{\fm}^{d-i}(M)^*
\end{equation}
as graded left $A$-modules
\cite[Theorem 4.2]{YZ}.

\subsection{The condition $\chi$}

Let $A$ be a connected graded $k$-algebra.

We say $A$ {\sf satisfies condition} $\chi$
if $\uExt^i_A(k,M)$ is finite dimensional for all finitely
generated $M$ and all $i$.
By \cite[Cor. 3.6]{AZ}, this is equivalent to $H^i_{\fm}(M)$ being
zero in large positive degree for all $i$ and all finitely generated 
$M$. 
Hence, if $A$ is noetherian and Cohen-Macaulay,
formula (\ref{dual}) implies that $A$ satisfies $\chi$ on both sides
\cite[Theorem 4.2]{YZ}.
The precise relationship between condition $\chi$ and the
Cohen-Macaulay property is given by \cite[Theorem 6.3]{vdB7}.

A noetherian, connected, graded algebra $A$ satisfying 
$\chi$ has finite depth, and for every $M
\in \mod  A$ of finite projective dimension, 
$$
\pdim M + \depth M = \depth A.
$$
As in the commutative case we call this the Auslander-Buchsbaum
formula. The non-commutative version was proved by J\"orgensen
\cite{J1}.

\subsection{Non-commutative spaces}

Let $A$ be a connected graded noetherian $k$-algebra. 

Artin and Zhang  \cite{AZ} define $\Proj A$ to be the (imaginary) non-commutative scheme defined implicitly
by declaring that the category of ``quasi-coherent sheaves'' on it is
$$
 \Qcoh \big( \Proj A\big):=\QGr A:= \frac{\Gr A}{\Fdim A}
$$
where $\Fdim A$ is the full subcategory 
consisting of direct limits of finite dimensional modules.
We write $\pi:\Gr A \to \Qcoh (\Proj A)$ for the quotient functor 
and $\omega$ for its right adjoint. Modules in $\Fdim A$ are said to be 
{\sf torsion}. 

We also define
$$
\coh  \big( \Proj A\big):= \qgr A:= \frac{\gr A}{\Fdim A \cap \gr A}.
$$
It is the full subcategory of noetherian objects in $ \Qcoh \big( \Proj A\big)$.

Write $X=\Proj A$ and $\cO_X=\pi A$.

Artin and Zhang define the cohomology groups
$$
H^q(X,\cF):=\Ext^q_X(\cO_X,\cF).
$$
If $M$ is a graded $A$-module, there is an exact sequence
\begin{equation}
\label{eq.omega.pi}
0 \to H^0_{\fm}(M) \to M \to \omega\pi M \to H^1_{\fm}(M) \to 0,
\end{equation}
and, if $\cM=\pi M$, then
\begin{equation}
\label{eq.HH}
H^q(X,\cM) \cong H^{q+1}_{\fm}(M)_0
\end{equation}
for $q \ge 1$ \cite[Prop. 7.2]{AZ}. 

Following \cite[Defn. 2.4]{YZ}, we say that $X$ is {\sf
Cohen-Macaulay} of dimension $d$ if there exists
$\omega_X \in \coh X$ and isomorphisms
$$
H^q(X,-) \to \Ext^{d-q}_X(-,\omega_X)^*
$$
on $\coh X$ for all $q$.

Suppose  $A$ is noetherian and Cohen-Macaulay of depth $d+1$.
Since $A$ satisfies $\chi$, (\ref{eq.HH}) allows us to quote
\cite[Thm. 2.3]{YZ} which says that $X=\Proj A$ is Cohen-Macaulay 
of dimension $d$ with $\omega_X\cong \pi(H_{\fm}^{d+1}(A)^*)$.

\subsection{Non-commutative analogues of $\PP^n$}
\label{sect.qu.Pn}

The quadric surfaces of interest to us are degree two hypersurfaces
in quantum $\PP^3$s.

For the purposes of this paper a {\sf quantum $\PP^n$} is a non-commutative scheme $\Proj S$ 
for which $S$ is a connected graded
$k$-algebra with the following properties:
\begin{enumerate}
\item{}
$S$ has global homological dimension $n+1$ on both sides and
$$
\Ext^i(k,S) \cong 
	\begin{cases} 
	0 & \text{if $i \ne n+1$}
	\\
	k & \text{if $i = n+1$}
	\end{cases}
$$
for the right and left trivial modules $k=S/S_{\ge 1}$
(i.e., $S$ is an {\sf Artin-Schelter (AS) regular} algebra);
\item{}
$S$ is right and left noetherian;
\item{}
$H_S(t)=(1-t)^{-n-1}$.
\end{enumerate}
J.J. Zhang showed that these conditions imply that $S$ is a
Koszul algebra and has dualizing module 
$\omega_S \cong A(-n-1)$ \cite[Thm. 5.11]{Sm4}.
Furthermore, $S$ satisfies $\chi$ on both sides. When $n+1 \le 4$,
$S$ is a domain by \cite{ATV}.

A result of Shelton and Vancliff \cite[Lemma 1.3]{ShV} shows that
for quantum $\PP^3$s the hypotheses (1)-(3) are not the most
efficient---one can slightly weaken them.

Write $\PP^n_{nc}=\Proj S$. The hypotheses ensure that
$H^{n+1}(\PP^n_{nc},-)=0$ and that the dimensions of $H^q(\PP^n_{nc},
\cO(r))$ agree with those in the commutative case.

The Grothendieck group of a quantum $\PP^n$ is isomorphic
to $\ZZ[t,t^{-1}]/(1-t)^{n+1}$ with $[\cF(-1)]=[\cF]t$.
There is a good notion of degree for closed subspaces of $\Proj
S$. In particular, if $z \in S$ is a homogeneous normal element,
meaning that $Sz=zS$, then $\Proj S/(z)$ is a hypersurface of
degree equal to $\deg z$. 
Write $A=S/(z)$. Then $A$ is Gorenstein of depth $n$, and satisfies
$\chi$. In particular, $H^n(\Proj A,-)=0$.

\begin{lemma}
If $S$ is a connected graded $k$-algebra of finite global dimension
and $H_S(t)=(1-t)^{-n}$, then the Hilbert series of every 
finitely generated $A$-module of GK-dimension one is
eventually constant.
\end{lemma}
\begin{pf}
The minimal projective resolution of a finitely generated $A$-module
$M$ is finite, and all terms are direct sums of shifts of $A$, so the 
Hilbert series of the module is of the form $f(t)(1-t)^{-n}$ for 
some $f(t) \in \ZZ[t,t^{-1}]$. The hypothesis on the GK-dimension
means that we can rewrite this as $g(t)(1-t)^{-1}$ with 
$g(t) \in \ZZ[t,t^{-1}]$. Hence $\dim M_n=g(1)$ for $n \gg 0$.
\end{pf}

\section{Maximal Cohen-Macaulay modules}

Suppose  $A$ is a connected, graded,  noetherian, and Cohen-Macaulay of depth $d
\ge 1$.

\subsection{}
A noetherian $A$-module $M$ is {\sf maximal Cohen-Macaulay}
if $\depth M=d$. We write $\MCM(A)$ for the full subcategory of 
$\gr A$ consisting of the maximal Cohen-Macaulay modules; we consider
the zero module to be maximal Cohen-Macaulay,
so $\MCM(A)$ is an additive category.

If $i \ge d$, then $\Omega^iM \in \MCM(A)$ for all $M \in \gr A$. If $M$ is in $\MCM(A)$ so
is $\Omega M$.  

The {\sf  stable category of maximal Cohen-Macaulay modules}, denoted $\uMCM(A)$, has the same
objects as $ \MCM(A)$ and morphisms 
$$
\Hom_{\uMCM(A)}(M,N):= \frac{ \Hom_{\Gr A}(M,N)}{P(M,N)}
$$
where $P(M,N)$ consists of the degree zero $A$-module 
maps $f:M \to N$ that factor through a projective in 
$\Gr A$. 
%Projective modules become isomorphic to zero in
%$\uMCM(A)$ because the identity morphism on a projective belongs to
%$P$. Hence, if $M$ is a finitely generated
%module of depth $n$ and $i+n \ge d$, the $i^{\th}$ syzygy of 
%$M$ is a well-defined object of $\uMCM(A)$ up to isomorphism.

\subsection{}
The next two results are due to Buchweitz and are stated in his appendix to 
the paper \cite{BEH}; see also \cite{Buch}. 

\begin{theorem}
Suppose $A$ is Gorenstein. Then
$\uMCM(A)$ is a triangulated category with respect to the 
translation functor $M[-1]:=\Omega M$.
If $M$ and $N$ are maximal Cohen-Macaulay modules, then 
$$
\Hom_{\uMCM(A)}(M,N[n]) \cong \Ext^n_{\Gr A}(M,N)
$$
for all $n \ge 1$.
\end{theorem}

\begin{theorem}
\label{thm.Buch}
Let $A$ be a Gorenstein, connected, graded, Koszul algebra over 
a field $k$, and $A^!$ its quadratic dual. 
The Koszul duality functor $K$ fits into a commutative diagram
\begin{equation}
\label{eq.MCM.proj}
\xymatrix @=3pc 
{
\MCM(A)\ar[r] \ar[d] & \gr A \ar[r] & \sD^b(A) \ar^{K}[r] &
\sD^b(A^!)  \ar[d]
\\
\uMCM(A) \ar[rrr]_B &&& \sD^b(\qgr A^!)
}
\end{equation}
in which the bottom arrow is a duality
$$
\uMCM(A)  \cong \sD^b(\qgr A^!),
\qquad M \mapsto \RuHom_A(M,k).
$$
The $t$-structure on $\uMCM(A)$ induced by 
the natural $t$-structure on $\sD^b(\qgr A^!)$  is
\begin{align*}
\uMCM(A)^{\ge p} &=\{M \; | \; \uExt^i_A(M,k)_j=0 \hbox{ for } i+j >p\}
\\
\uMCM(A)^{\le p} &=\{M \; | \; \uExt^i_A(M,k)_j=0 \hbox{ for } i+j <p\}
\end{align*}
The heart for this $t$-structure consists of the maximal Cohen-Macaulay modules having a
linear resolution.
\end{theorem}

We will refer to the duality
$$
B:\uMCM(A) \to \sD^b(\qgr A^!)
$$
in Theorem \ref{thm.Buch} as ``Buchweitz's duality''.

\begin{lemma}
\label{lem.F.twists}
Suppose $A$ is a connected graded, Gorenstein, Koszul algebra.
Write $F$ for the composition
\begin{equation}
\begin{CD}
\gr A @>>> \sD^b(A) @>{K}>> \sD^b(A^!) @>>> \sD^b(\qgr A^!)
\end{CD}
\end{equation}
If $M \in \gr A$, then $F(\Omega M) \cong (FM)[1]$ and
$F(\Omega M(1)) \cong (FM)(1)$.
\end{lemma}
\begin{pf}
There is an exact sequence $0 \to \Omega M \to \oplus_{i \in I}
A(i) \to M \to 0$ for some multiset $I$, and hence a distinguished
triangle $\Omega M \to \oplus_{i \in I}  A(i) \to M \to $ in 
$\sD^b(A)$. 
The image of $K(\oplus A(i))$ in $\sD^b(\qgr A^!)$ is zero so
there is an isomorphism
$$
\begin{CD}
\pi KM @<{\sim}<< \pi K((\Omega M)[1]) \cong \pi K(\Omega M)[-1]
\end{CD}
$$
in $\sD^b(\qgr A^!)$.
Hence $FM \cong F(\Omega M)[-1]$.

The other isomorphism is established in a similar way.
\end{pf}

\begin{lemma}
\label{lem.indec.MCM}
Suppose $A$ is a connected graded, Gorenstein algebra.
The isomorphism classes of indecomposable objects in $\uMCM(A)$ are
in bijection with the isomorphism classes of indecomposable 
non-projective modules in $\MCM(A)$.
\end{lemma}
\begin{pf}
Non-trivial direct summands of an object in an additive category
correspond to non-trivial idempotents in its endomorphism ring.

Let $M$ be an indecomposable non-projective in $\MCM(A)$.
Then $E=\End_{\Gr A}M$ is a finite dimensional local ring, meaning
that $E/\rad E$ is a division ring.
Hence the endomorphism ring of $M$ in $\uMCM(A)$ is also local, so $M$
is indecomposable in $\uMCM(A)$.

Let $M'$ be another indecomposable non-projective in $\MCM(A)$, and
suppose that $f:M \to M'$ and $g:M' \to M$ become mutually inverse
isomorphisms in $\uMCM(A)$. To show that $M$ is isomorphic to $M'$ in
$\MCM(A)$, it suffices to show that $fg$ and $gf$ are isomorphisms 
in $\MCM(A)$. It therefore suffices to show that if
$h:M \to M$ is an isomorphism in $\uMCM(A)$, then it is an isomorphism
in $\MCM(A)$. But this is clear, since the isomorphisms $M \to M$ in
either category are the endomorphisms that are not in the
radical. 

We have shown that the functor $\MCM(A) \to \uMCM(A)$
gives an injective map from the set of isomorphism classes of
indecomposable non-projective Cohen-Macaulay modules to the set of
isomorphism classes of indecomposable objects in $\uMCM(A)$. 
We now show this map is surjective. 
If $M$ is a maximal Cohen-Macaulay module
that becomes indecomposable as an object in $\uMCM(A)$ we may 
write $M$ as a direct sum of indecomposables in $\MCM(A)$
and this gives a direct sum decomposition of $M$ in $\uMCM(A)$ each
term of which is either zero or indecomposable; hence, in $\uMCM(A)$,
$M$ is isomorphic to some $M'$ where $M'$ is an indecomposable
non-projective in $\MCM(A)$.
\end{pf}

{\bf Remark.}
Suppose $A$ is Gorenstein, and let $N$ be a maximal 
Cohen-Macaulay module having no non-zero projective direct summand. 
By applying $\uHom_A(-,A)$ to $0 \to \Omega N \to P \to N \to 0$,
where $P \to N \to 0$ is the start of a minimal projective resolution,
one sees that $\Omega N$ also has no non-zero  projective direct
summand. Hence for all $d \ge 0$, $\Omega^d N$ has no
non-zero  projective direct summands.

\subsection{Simple objects and maximal Cohen-Macaulay modules}

\begin{lemma}
\label{lem.simple.MCM}
Let $A$ be a connected, graded, Gorenstein, Koszul algebra and 
$$
B:\uMCM(A) \to \sD^b(\qgr A^!)
$$
 the equivalence in Theorem 
\ref{thm.Buch}. Let $S$ be a simple object in $\qgr A^!$. Then 
\begin{enumerate}
\item{}
there is a unique-up-to-isomorphism indecomposable $M \in \MCM(A)$ 
such that $BM \cong S[0]$; 
\item{}
if $S \cong S(d)$, then $\Omega^dM(d) \cong M$;
\item{}
$\Omega^nM(n)$ has a linear resolution for all $n$.
\end{enumerate}
\end{lemma}
\begin{pf}
The equivalence of categories and Lemma \ref{lem.indec.MCM} ensure
the existence and uniqueness of $M$. Because $S[0]$ is in the heart
of $\sD^b(\qgr A^!)$, $M$ has a linear resolution. An induction
argument using Lemma \ref{lem.F.twists} shows that $B(\Omega^dM(d))
\cong (BM)(d) \cong S(d) \cong S$, so $M \cong \Omega^dM(d)$ in 
$\uMCM(A)$. Because $M$ is indecomposable, it follows that $M$ is
isomorphic to a direct summand of $\Omega^dM(d)$ in $\MCM(A)$ and
the complementary summand is projective. By the previous remark, 
$\Omega^d M$ has no non-zero projective direct summand, so 
$\Omega^dM(d) \cong M$.

Since $M$ has a linear resolution, so does each $\Omega^n M(n)$.
\end{pf}

\section{Smoothness}

A non-commutative space 
$X$ is {\sf smooth} of global dimension $d$ if $d$ is the 
largest integer such that $\Ext_X^d(M,N) \ne 0$ for some $X$-modules
$M$ and $N$.

We now consider the question of what homological properties
of a connected graded $k$-algebra $A$ imply that $\Proj A$ is smooth.

\subsection{}
The following summarizes the commutative case.

\begin{proposition}
Let $A$ be a graded quotient of a positively graded  polynomial ring.
Write $X=\Proj A$.
The following are equivalent:
\begin{enumerate}
\item
$\gldim X \le d$;
\item
$\dim_k \uExt^i_A(M,N) < \infty$ for all $i \ge d+1$ and all $M,N \in
\gr A$;
\item
whenever $N \to E^\bullet$ is a minimal injective resolution in
$\Gr A$, $E^i$ is torsion for all $i \ge d+1$.
\end{enumerate}
\end{proposition}
\begin{pf}
(3) $\Leftrightarrow$ (1)
If  $0 \to N \to E^\bullet$ is a minimal injective resolution in
$\Gr A$, then $0 \to \pi N \to \pi  E^\bullet$ is a minimal
injective resolution in $\Qcoh X$. 
Thus $\gldim X  \le d$ if and only if $\pi E^i =0$ for all $i \ge
d+1$.
Hence the equivalence of conditions (1) and (3)

(2) $\Rightarrow$ (1)
Since $A$ is commutative, $\uExt^i_A(M,N)$ is an $A$-module. 
If $f$ is a homogeneous element of $A$ lying in $\fm$, then 
$A[f^{-1}] \otimes_A \uExt^i_A(M,N) =0$ for $i>d$, and so
$\uExt^i_{A[f^{-1}]} \equiv 0$ for $i>d$. Hence $\Proj A[f^{-1}]$
is smooth of global dimension at most $d$.
Since $X$ is covered by open affines of the form $\Spec
A[f^{-1}]_0 \cong \Proj A[f^{-1}]$ with $f \in \fm$, it follows
that $\gldim X \le d$.

(3) $\Rightarrow$ (2)
If condition (3) holds then applying $\uHom_A(M,-)$ to a minimal
injective resolution of $N$ produces a complex consisting of torsion
modules after the $d^{\th}$ term. Hence  $\uExt_A^i(M,N)$ is
torsion for $i>d$. However, if $M$ and $N$ are noetherian, then
$\uExt_A^i(M,N)$ is a noetherian $A$-module as one sees by applying
$\uHom_A(-,N)$ to a minimal projective resolution of $M$. Hence
$\uExt_A^i(M,N)$ is finite dimensional for $i>d$ whenever  $M,N
\in \gr A$.
\end{pf}

The proof of (3) $\Leftrightarrow$ (1) works when $A$ is
not commutative, but the other two parts of the proof fail because
$\uExt^i_A(M,N)$ is not an $A$-module when $A$ is not
commutative. Nevertheless, we will show that the implication (3) 
$\Rightarrow$ (2) holds if $A$ satisfies $\chi$.

First we need the following lemma that we learned from Kontsevich.

\begin{lemma}
\label{lem.gabber}
Let $A$ be a noetherian  connected graded $k$-algebra satisfying
$\chi$.
%Write $X=\Proj A$.
If $\gldim\big( \Proj A\big)=d<\infty$, and $M \in \gr A$, then there is
a perfect complex $P \in \sD^b(\gr A)$ concentrated in homological
degree $[-d,0]$ and a 
bounded complex $N$ together with a map in
$M_{\ge n} \oplus N \to  P$
($n$ large) whose cone has finite dimensional cohomology. 
\end{lemma}
\begin{pf}
Take an exact sequence $0 \to Z \to P_{d} \to \cdots \to P_0 \to
M \to 0$ with each $P_i$ a finitely generated free module.
Applying $\pi$ to this gives an element of $\Ext_X^{d+1}(\pi M,\pi Z)$
which must be zero, so the triangle  $\pi Z[d]\to \pi P \to 
\pi M \to $ is split. Applying $\omega$ gives an isomorphism 
$\omega\pi M \oplus  \omega\pi
Z[d] \to \omega\pi P$. However, since $\chi$ holds, for every
finitely generated module $N$ the map $N_{\ge
n} \to (\omega\pi N)_{\ge n}$ has finite dimensional cokernel (and
it obviously has finite dimensional kernel). Hence, for large $n$
there is an isomorphism $M_{\ge n} \oplus Z[d]_{\ge n} \cong P_{\ge n}$,
and hence a map $M_{\ge n} \oplus Z[d]_{\ge n} \to P$ whose cone has
finite dimensional cohomology.
\end{pf}

\begin{proposition}
Let $A$ be  a noetherian  connected graded $k$-algebra
satisfying $\chi$. 
Write $X=\Proj A$.
If $\gldim X=d<\infty$, then
$\dim_k\uExt_A^i(M,N)< \infty$ for all $i \ge d+1$ and all 
$M,N \in \gr A$.
\end{proposition}
\begin{pf}
By Lemma \ref{lem.gabber}, there is a distinguished triangle $M_{\ge
n} \oplus N \to  P \to C \to $ such that $\oplus_qH^q(C)$ is finite
dimensional.
\end{pf}

\subsection{}

We call a triangulated category is {\sf semisimple} if 
in every distinguished triangle 
$$
\begin{CD}
L @>u>> M @>v>> N @>w>>
\end{CD}
$$
at least one of $u$, $v$, and $w$, is zero. 
The condition that $w=0$ is equivalent to the condition 
that $u$ and $v$ induce an isomorphism $M \cong L \oplus N$.
It follows that if $L$ and $M$ are objects in a 
semisimple triangulated category, then 
$\Hom(L,M) \ne 0$ if and only if $\Hom(M,L) \ne 0$. 
Futhermore, the heart of every $t$-structure on a semisimple
triangulated category is semisimple, meaning that every short exact
sequence splits. The derived category of 
a semisimple abelian category is semisimple.

\smallskip

{\bf Notation.}
We write $\MCM_{\ge n}$ for the full subcategory of
$\MCM(A)$ consisting  of graded
Cohen-Macaulay $A$-modules $M$ such that $M_i=0$ for all $i<n$.
We write $\uMCM_{\ge n}$ for the essential image of $\MCM_{\ge n}$
in $\uMCM(A)$.

\begin{proposition}
\label{prop.mcm.smooth}
Let $A$ be a connected, graded $k$-algebra that is
Gorenstein and satisfies $\chi$.
If $\uMCM(A)$ is semisimple, then $\Proj A$ is smooth.
\end{proposition}
\begin{pf}
Suppose  $A$ has depth $d$, so the $d^{\th}$ syzygy of a
finitely generated module is maximal Cohen-Macaulay.

Fix $K \in \MCM(A)$. 
For all integers $n$ greater than the degrees of the 
minimal generators of $K$, $\Hom_{\uMCM}(K,-)$ vanishes on
$\uMCM_{\ge n}$. Since $\uMCM(A)$ is semisimple, it follows that
$\Hom_{\uMCM}(-,K)$ also vanishes on $\uMCM_{\ge n}$.

Now fix $M,N \in \gr A$.
Then
\begin{align*}
\Ext_{\Proj A}^{d+1}(\pi M,\pi N) & \cong \Ext^{d+1}_{\Gr A}(M_{\ge n},N)
\qquad \hbox{for $n \gg 0$} 
\\
& \cong \Ext^{1}_{\Gr A}(\Omega^d(M_{\ge n}),N) 
\\
& \cong \Ext^{d+1}_{\Gr A}( \Omega^d(M_{\ge n}),\Omega^dN)
\qquad \hbox{because $A$ is Gorenstein}
\\
& \cong \Hom_{\uMCM}( \Omega^d(M_{\ge n}),(\Omega^dN)[d+1]).
\end{align*}
The previous paragraph shows that this is zero for $n \gg 0$
because $\Omega^d(M_{\ge n}) \in \MCM_{\ge n}$.
\end{pf}

Part of the argument in Proposition \ref{prop.mcm.smooth}
can be restated in the following way.

\begin{proposition}
\label{prop.mcm.ss}
Suppose $\uMCM(A)$ is semisimple.
If $N \in \MCM(A)$, there is an integer $n$ such that
$\Ext^1_{\Gr A}(M,N) =0$ for all $M \in \MCM_{\ge n}$.
\end{proposition}
\begin{pf}
Choose an integer
$
n>\{\hbox{the degrees of a minimal set of generators for $N$}\}.
$
Suppose  $M \in \MCM_{\ge n}$. Then $\Omega M \in \MCM_{\ge n}$
also. Hence $\Hom_{\Gr A}(N,\Omega M)$ is zero, and so is its
quotient $\Hom_{\uMCM(A)}(N,M[-1])$. But $\uMCM(A)$ is
semisimple, so $\Hom_{\uMCM(A)}(M[-1],N)=0$ also. 
Thus $\Ext^1_{\Gr A}(M,N) =0$.
\end{pf}

\begin{proposition}
\label{prop.mcm.hyper}
Let $S$ be a Gorenstein $k$-algebra of finite global 
dimension and $z$ a central regular element 
of degree $d$. Let $A=S/(z)$. If $M \in \MCM(A)$ is not projective, then
\begin{enumerate}
\item{}
there is a resolution $ 0 \to S^s \to S^s \to M \to 0$ of ungraded
$S$-modules;
\item{}
$\Omega^2M \cong M(-d)$;
\item{}
$ \uExt^{i+2}_A(M,N) \cong \uExt_A^i(M,N)(d) $
for all $A$-modules $N$ and all $i \ge 1$.
\end{enumerate}
\end{proposition}
\begin{pf}
(1)
We have $\depth_S M=\depth_AM =\depth_A A=\depth_S A = \depth S -1$,
so $\pdim_S M=1$ by the Auslander-Buchsbaum formula. 
Hence $M$ has a free resolution  $ 0 \to S^r \to S^s \to M \to 0$.
But $r \le s$ because $S$ is noetherian and $s \le r$ because $S^sz
\subset S^r$, so $s=r$.
%and $r=s$ because $M$ is an $A$-module.

Although this is a resolution of $M$ in $\QGr S$, the minimal
resolution of $M$ in $\Gr S$ also has this form although we may
have to change $s$ and we will need to place some gradings on the
two free modules.

(2) and (3)
From the presentation of $A$ we get  $\Tor^S_1(M,A) \cong M(-d)$. 
Now, by applying $-\otimes_S A$ to the resolution of $M_S$, we
obtain an exact sequence
$$
0 \to M(-d) \to A^s \to A^s \to M \to 0
$$
in $\Mod A$.
We can give the two copies of $A^s$ gradings so this becomes a
sequence of graded $A$-modules. Thus $\Omega^2M \cong M(-d)$.
The result now follows by dimension-shifting.
\end{pf}

\section{Quadrics in quantum $\PP^3$s}

\subsection{The algebra $C(A)$}
Using the notation in part (2) of the next lemma, we define
\begin{equation}
\label{defn.C(A)}
C(A):=A^![w^{-1}]_0.
\end{equation}

We write $\Mod \, C$ for the category of right modules over a ring $C$, 
and $\mod \, C$ for the full subcategory of finitely presented modules.

\begin{lemma}
\label{lem.C(A)}
Let $S$ be a connected, graded, noetherian, Koszul algebra 
of finite global dimension, $z$ a central regular element 
of degree two, and $A=S/(z)$. 
Then
\begin{enumerate}
\item{}
$A$ is a Koszul algebra;
\item{}
there is a central, regular element $w \in A_2^!$ such that 
$A^!/(w)=S^!$;
\item{}
the algebra $C(A)$ has finite dimension equal to 
$\dim_k (S^!)^{(2)}$, the dimension of the even degree
part of $S$;
\item{}
the categories $\QGr A^!$ and $\Mod C(A)$ are equivalent
via $\pi N \mapsto N[w^{-1}]_0$, where $N \in \Gr A^!$.
\end{enumerate}
\end{lemma}
\begin{pf}
(1) and (2).
The proof is similar to that for modding out a central regular element
of degree one \cite{LSV}.

(3)
Because $S$ is Koszul, the hypothesis that
$\gldim S < \infty$ implies that $S^!$ has finite dimension.
It follows that $A^!_{m+2}=wA^!_m$ for large $m$, and hence that 
$$
A^![w^{-1}]_0= A^!_0+A^!_2w^{-1}+ \cdots =A_{2n}^!w^{-n}
$$
 for $n \gg 0$. In particular,
$\dim_k A[w^{-1}]_0=\dim_k A_{2n}$ for $n \gg 0$.
By (1) and (2),  
$(1+t)H_{A^!}(t)=(1-t)^{-1}H_{S^!}(t)$ so
$$
\dim_k A^!_{2n}=\dim_k S^!_0+ \dim_k S^!_2+\cdots,
$$
for $n \gg 0$, as required.

(4)
A graded $A^!$-module has finite dimension 
if and only if it is annihilated by a power of $w$, 
so $\Gr A^!/\Fdim A^!$ is equivalent to 
$\Gr A^![w^{-1}]$. Since $A^!$ is generated in degree one, 
$A^![w^{-1}]$
is strongly graded, and therefore 
$\Gr A^![w^{-1}]$ is equivalent to $\Mod A^![w^{-1}]_0$.
\end{pf}

The degree shift functor $(1)$ on $\Gr A^!$ induces
auto-equivalences of $\QGr A^!$ 
and $\Mod C(A)$ that we still denote by $(1)$.
Since $w$ is central and homogeneous of degree two, on $\Mod C(A)$ we have
$(2) \cong \id_{\Mod C(A)}$. 

Notice that $A^!$ is noetherian because $A^!/(w)$ is.

\begin{proposition}
\label{prop.C(A).smooth}
Let $S$ be a Gorenstein, connected, graded, noetherian, Koszul algebra 
of finite global dimension, $z$ a central regular element 
of degree two, and $A=S/(z)$. 
\begin{enumerate}
\item{}
There are equivalences of categories
$$
\setlength{\unitlength}{1mm}
\begin{picture}(25,15)
\put(-11,13){$\uMCM(A)$}
\put(4,14){\vector(1,0){20}}
\put(14,15){$B$}
\put(25,13){$\sD^b(\qgr A^!)$} 
\put(3,11){\vector(1,-1){7}}
\put(24,11){\vector(-1,-1){7}}
\put(5,0){$\sD^b(\mod C(A))$}
\end{picture}
$$
\item{}
If $C(A)$ is a semisimple ring, then $\Proj A$ is smooth.
\end{enumerate}
\end{proposition}
\begin{pf}
The horizontal equivalence is given by Theorem \ref{thm.Buch},
and the southwest equivalence is given by Lemma \ref{lem.C(A)}.
Part (2) follows from (1) and Proposition \ref{prop.mcm.smooth}
because the derived category of a semisimple abelian category is 
semsimple, hence abelian.
\end{pf}

\subsection{Notation and Hypotheses}
\label{sect.notn.hyp}
We fix the following hypotheses and notation for the remainder of this section:
$k$ is an algebraically closed field,
$S$ denotes a connected, graded, noetherian, Gorenstein, Koszul algebra with Hilbert series
$(1-t)^{-4}$; $z$ is a non-zero, homogeneous, central element of 
degree two such that $A:=S/(z)$ is a domain. 
We write $Q:=\Proj A$.

The hypotheses imply that $H_{S^!}(t)=(1+t)^4$, so $\gldim S=4$.
The previous two results apply, so the finite dimensional algebra $C(A)$
is well-defined.

It follows from Corollary \ref{cor.ausl} below and \cite[Thm.
3.9]{ATV} that $S$ is a domain.
Thus $\Proj S$ is a quantum $\PP^3$ and $Q$ is a 
quadric hypersurface in it. The assumption that $A$ is a domain says that
$Q$ is ``reduced and irreducible''.

\begin{proposition}
\label{prop.C(A)}
Suppose  $S$ is a connected, graded, noetherian, Gorenstein, Koszul algebra  such that
$H_S(t)=(1-t)^{-4}$. Let $0 \ne z \in S_2$ be a central element and
suppose that $A:=S/(z)$ is a domain. Then 
\begin{enumerate}
\item{}
$\dim_k C(A)=8$;
\item{}
$C(A)$ has no one-dimensional modules;
\item{}
the following are equivalent:
\begin{enumerate}
\item{}
$C(A)$ is semisimple;
\item{}
$C(A)$  has two simple modules up to isomorphism;
\item{}
$C(A) \cong M_2(k) \oplus M_2(k)$.
\end{enumerate}
\end{enumerate}
\end{proposition}
\begin{pf}
(1)
This does not depend on $A$ being a domain.
Because $z$ is regular, 
$H_A(t)=(1+t)(1-t)^{-3}$ and $H_{A^!}(t)= (1-t)^{-1}(1+t)^3$.
Thus $\dim_k A_n=8$ for $n \gg 0$, and 
$
\dim_k C(A)=8
$
by Lemma \ref{lem.C(A)}.

(2)
First we show that if $V$ is a subspace of $A^!_1=A_1^*$ of
codimension one, then $A_1^!V=VA_1^!=A_2^!.$
To see this write $V=a^\perp$ where $a \in A_1$;
because $A$ is a domain, $a \otimes A_1 \cap R
=0$, where $R$ denotes the relations in $A_1 \otimes A_1$
defining $A$; hence $a^\perp \otimes A_1^* +R^\perp = 
A_1^* \otimes A^*_1$; thus $VA_1^!=A_2^!$. Similarly,
$A_1^!V=A_2^!$.

\underline{Claim:}
Let $T$ be a connected graded algebra generated in degree one, and
$w\in T_d$ a central regular element of degree $d>1$. 
If $T_2=T_1V$ for every
codimension one subspace $V \subset T_1$, then $T[w^{-1}]_0$ does
not have a one-dimensional module.

\underline{Proof:}
Suppose to the contrary that $N_0$ is a one-dimensional
$T[w^{-1}]_0$-module. Then $N_0$ is the degree zero component of
the $T[w^{-1}]$-module $N:=T[w^{-1}] \otimes_{T[w^{-1}]_0} N_0$.
Because $w$ is a unit, $N_{di}=w^iN_0$ for all $i \in \ZZ$. In
particular, $\dim_k N_{di}=1$ for all $i$.
Hence if $m \in N_{di-1}$, then $Vm=0$ for some subspace $V \subset
T_1$ of codimension at most one. Hence $T_2m=0$. It follows that
$$
0=T_2N_{di-1}=T_3N_{di-2}=\cdots = T_{d+1} N_{di-d}.
$$
In particular, $T_{d+1} N=0$, so $w^2N=0$, contradicting the fact
that $w$ is a unit in $T[w^{-1}]$. $\lozenge$

The claim applies to $T=A^!$, so (2) follows.

(3)
This follows from (1) and (2) because $k$ is algebraically closed.
\end{pf}

\subsection{The set $\MM$}
\label{sect.notn.MM}
We define
$$
\MM:=\{M \in \MCM(A) \; | \; \hbox{$M$ is indecomposable, $M_0
\cong k^2$, $M=M_0A$}\}.
$$

Because $A$ is a noetherian domain it has a division ring of
fractions, say $Q$, and we may define the {\sf rank} of an
$A$-module $N$ as $\dim_Q N \otimes_A Q$. 

\begin{proposition}
\label{prop.min.res.MCM}
Let $A$ and $S$ be as in Proposition \ref{prop.C(A)}.
If $M \in \MM$, then
\begin{enumerate}
\item{}
$M \cong \Omega^2M(2)$;
\item{}
its minimal resolution is 
$\cdots \to A(-2)^2 \to A(-1)^2 \to A^2 \to M \to 0$;
\item{}
$\rank M=1$;
\item{}
$H_M(t)=2(1-t)^{-3}$;
\item{} 
$M$ is 3-critical with respect to GK-dimension.
\end{enumerate}
Furthermore, there is a bijection
\begin{align*}
\MM & \longleftrightarrow \{\hbox{simple $C(A)$-modules}\}
\\
M & \longleftrightarrow FM.
\end{align*}
\end{proposition}
\begin{pf}
The hypotheses on $M$ ensure that it is not projective. 

(1)
This was already established in Proposition \ref{prop.mcm.hyper}.

(2)
By (1), the minimal resolution of $\Omega^2 M$ begins $A(-2)^2 \to
\Omega^2 M \to 0$. Combining this with Proposition
\ref{prop.mcm.hyper}, we see that the minimal resolution of $M$ begins
$$
\cdots \to A(-2)^2 \to A(-i) \oplus A(-j) \to A^2 \to M \to 0
$$
for some $i,j$. However, the minimality of the resolution
forces $i=j=1$. Since $\Omega^2 M \cong M(-2)$, the full minimal
resolution of $M$ can be constructed by splicing together shifts of
the exact sequence $0 \to M(-2) \to A(-1)^2 \to A^2 \to M \to 0$.

(3)
Because $A$ is a domain, its rank is one.
The result now follows from the exact sequence 
$0 \to M(-2) \to A(-1)^2 \to A^2 \to M \to 0$.

(4) 
This follows from (2).

(5)
Because $M(-1)$ embeds in $A^2$, every non-zero submodule of $M$
has GK-dimension three; since the rank of $M$ is one, all its
non-zero submodules have rank one too.
Hence every proper quotient of $M$ has GK-dimension $\le 2$.

We now establish the bijection between $\MM$ and the simple
$C(A)$-modules. Let $M \in \MM$. By (2), 
$M$ has a linear resolution, so $FM=N[0]$ for some $C$-module $N$. 
But $N=\uExt_A^*(M,k)[w^{-1}]_0
\cong \uExt_A^{2i}(M,k)$ for $i \gg 0$, so $\dim_k N=2$, and
Proposition \ref{prop.C(A)} now implies that $N$ is simple.

Conversely, let $N$ be a simple $C$-module. 
By Lemma \ref{lem.simple.MCM}, there
is a unique indecomposable maximal Cohen-Macaulay module $M$ such that
$FM \cong N[0]$, and $M$ has a linear resolution.
By  Proposition \ref{prop.C(A)}, $\dim_k N=2$. 
But $N \cong \uExt^*_A(M,k)[w^{-1}]_0$, so
$\uExt^{2i}_A(M,k) \cong k^2$ for $i \gg 0$. Hence $\Omega^{2i}M$
is generated by two elements for $i \gg 0$. But $\Omega^{2i}M \cong
M(-2i)$ by Proposition \ref{prop.mcm.hyper}, so $M$ is generated by
two elements, and these are of degree zero because $M$ has a linear
resolution. Hence $M \in \MM$.
\end{pf}

\begin{lemma}
\label{lem.two.mcms}
Let $A$ and $S$ be as in Proposition \ref{prop.C(A)}.
If $M \in \MM$, then
\begin{enumerate}
\item{}
$\Omega M(1) \in \MM$;
\item{}
there is an exact sequence $0 \to M(-1) \to A^2 \to \Omega M(1) \to 0$;
\item{} 
if $C(A)$ is not semisimple, then $M \cong \Omega M(1)$;
\item{}
if $M \not\cong \Omega M(1)$, then
$ \Hom_{\Gr A}(M,\Omega M(1))=0.$
\end{enumerate}
\end{lemma}
\begin{pf}
(1)
By the remark after Lemma \ref{lem.indec.MCM}, $\Omega M(1)$ is
indecomposable. From the minimal resolution for $M$ we see that
$\Omega M(1)$ is generated in degree zero and that $\dim_k \Omega
M(1)_0=2$.

(2)
This follows from the fact that $\Omega^2M \cong M(-2)$.

(3)
If $C(A)$ is not semisimple it has only one simple module so, by
(1) and the bijection in Proposition \ref{prop.min.res.MCM}, $M
\cong \Omega M(1)$.

(4)
A non-zero degree-zero homomophism $\a:M \to \Omega M(1)$ would be 
injective because $M$ and $\Omega M(1)$ are 3-critical with 
respect to GK-dimension, so its restriction $M_0 \to
(\Omega M(1))_0$ would be an isomorphism, whence $\a$ would be 
surjective.  This would contradict the hypothesis that 
$M \not\cong \Omega M(1)$.
\end{pf}

\subsection{Smoothness of $Q$ and semisimplicity of $C(A)$}

\begin{theorem}
\label{thm.main}
The noncommutative quadric $Q$ is smooth if and only if $C(A)$ is
semisimple.
\end{theorem}
\begin{pf}
($\Leftarrow$)
This was proved in Proposition \ref{prop.C(A).smooth}.

($\Rightarrow$)
We shall prove the contrapositive, so suppose that $C(A)$ is not
semisimple. Then it has only one simple module, $N$ say. Let $M(1) \in
\MM$ be such that $F(M(1)) \cong N[0]$. Thus $M$ is generated in
degree one.

We write $\cM$ for $\pi M$.

By Lemma \ref{lem.two.mcms}, $M \cong \Omega M(1)$, so
there is an exact sequence 
\begin{equation}
\label{ses.MAM(1)}
\begin{CD}
0 @>>> M @>>> A^2 @>>> M(1) @>>> 0.
\end{CD}
\end{equation}
This gives an exact sequence 
\begin{equation}
\label{ses.MAM(1).Q}
\begin{CD}
0 @>>> \cM @>>> \cO_Q^2 @>>> \cM(1) @>>> 0
\end{CD}
\end{equation}
in $\Mod Q$.

The result that $\gldim Q=\infty$ will be established in Step 3
below.

\underbar{Step 1.}
$\Ext^1_Q(\cM,\cM) \ne 0.$

\underbar{Proof:}
Because $M$ is maximal Cohen-Macaulay 
$$
\Ext^1_{\Gr A}(M,M) \cong \Hom_{\uMCM(A)}(M,M[1]).
$$
Applying the contravariant equivalence $F$ this is isomorphic
to
$$
\Hom_{\sD(C(A))}(F(M[1]),FM) \cong \Hom_{\sD(C(A))}(FM,(FM)[1]).
$$
But $FM$ is a translate of the unique simple $C$-module $N$, so
%$\Hom_{\sD(C(A))}(F(M),FM[1]) \cong 
the last term is isomorphic to $\Ext^1_{C(A)}(N,N)$ which
is non-zero because $C(A)$ is not semisimple.
Hence $\Ext^1_{\Gr A}(M,M) \ne 0$.

Applying $\pi$ to a non-split exact sequence 
\begin{equation}
\label{ses.MDM}
0 \to M \to D \to M \to 0
\end{equation}
in $\Gr A$ gives an exact sequence 
\begin{equation}
\label{ses.MDM.Q}
0 \to \cM \to \cD \to \cM \to 0
\end{equation}
in $\Mod Q$. 
Write $\theta:\MCM(A) \to \Gr A$ for the inclusion functor. 
There is an exact sequence of functors $0 \to H^0_{\fm} \to
\id_{\Gr A} \to \omega\pi \to H_{\fm}^1 \to 0$. These two
local cohomology functors vanish on $\MCM(A)$, so there is an
isomorphism of functors $\theta \to \omega\pi \theta$.
Hence, if (\ref{ses.MDM.Q}) were to split via a map $g:\cM \to \cD$, 
then $\omega(g)$ would provide a splitting of (\ref{ses.MDM}). It
follows that  $\Ext^1_Q(\cM,\cM) \ne 0.$

\underbar{Step 2.}
$\Ext^2_Q(\cM(1),\cM) \ne 0.$

\underbar{Proof:}
Applying $\Hom_{Q}(-,\cM)$ to (\ref{ses.MAM(1).Q}) gives an exact
sequence
$$
\Ext^1_Q(\cO_Q,\cM)^2 \to \Ext^1_Q(\cM,\cM) \to
\Ext^2_Q(\cM(1),\cM).
$$
The first term is zero because $\depth M=3$ implies that
$$
0= H^2_{\fm}(M) \cong H^1(Q,\cM) = \Ext^1_Q(\cO_Q,\cM).
$$
The second term is non-zero by Step 1, so the third term is
also non-zero, as required.

\underbar{Step 3.}
$\Ext^{n}_Q(\cM(n-1),\cM) \ne 0$ for all $n \ge 2$.

\underbar{Proof:}
We argue by induction on $n$. The case $n=2$ has already been
established in Step 2.
Applying $\Hom_{Q}(\cM(n),-)$ to (\ref{ses.MAM(1).Q}) gives an exact
sequence
\begin{equation}
\label{indn.seq}
\Ext_Q^n(\cM(n),\cO_Q)^2 \to \Ext_Q^n(\cM(n),\cM(1)) \to
\Ext_Q^{n+1}(\cM(n),\cM)
\end{equation}
The first term is isomorphic to two copies of
$\Ext_Q^n(\cM(n-2),\cO_Q(-2))$ which is isomorphic to
$H^{2-n}(Q,\cM(n-2))^*$ by Serre duality. This is zero for $n \ge
3$, and if $n=2$ it is isomorphic to $\Hom_{\Gr A}(A,M)^*=M_0$
which is zero because $M$ is generated in degree one. Since the
first term of (\ref{indn.seq}) is zero for all $n \ge 2$, we see from
the other two terms that the induction argument goes through.
\end{pf}

\begin{corollary}
\label{cor.M.OmegaM}
$C(A)$ is semisimple if and only if $M \not\cong \Omega M(1)$ for
all $M \in \MM$.
\end{corollary}
\begin{pf}
($\Leftarrow$)
If $C$ were not semisimple it would have a unique simple module
so, up to isomorphism, there would be only one module in
$\MM$; but if $M$ is in $\MM$ so is $\Omega M(1)$, whence $M \cong
\Omega M(1)$.

($\Rightarrow$)
If $C$ is semisimple, then $\gldim Q < \infty$. But the proof of
Theorem \ref{thm.main} showed that if there were an $M$ in $\MM$
such that $M \cong \Omega M(1)$, then $\gldim Q=\infty$. Hence
there can be no such $M$.
\end{pf}

\begin{corollary}
\label{cor.M.OmegaM.2}
If $Q$ is smooth, then $\MM$ consists of two non-isomorphic modules, 
say $\MM=\{M,M'\}$, and there are exact sequences
$$
0 \to M(-1) \to A^2 \to M' \to 0
$$
and
$$
0 \to M'(-1) \to A^2 \to M \to 0.
$$ 
\end{corollary}
\begin{pf}
This follows immediately from Lemma \ref{lem.two.mcms} and Corollary
\ref{cor.M.OmegaM}.
\end{pf}

\section{The Auslander property}

We fix the following notation in this section:
$S$ denotes a connected, graded, Gorenstein, Koszul algebra with Hilbert series
$(1-t)^{-4}$; $z$ is a non-zero, homogeneous, central element of 
degree two and $A:=S/(z)$.

The main result in this section, Theorem \ref{thm.ausl},
shows that $A$ has the {\sf Auslander property} by which we
mean that if $M \in \gr A$ and $N$ is a graded submodule of
$\uExt^j_A(M,A)$ for some $j$, then $\uExt^i_A(N,A)=0$ for $i<j$.
By \cite{Lev}, this will imply that $S$ also has the Auslander
property.

\subsection{}
The {\sf grade} of $M \in \gr A$ is $j(M):=\inf\{j \; | \;
\uExt^j_A(M,A) \ne 0\}$.
The Auslander property is equivalent to the condition that $j(N)
\ge j$ for all submodules $N \subset \uExt^j_A(M,A)$.
To prove that $A$ has the Auslander property we first
prove that
$$
j(M)+\GKdim M =3
$$
for all $M \in \gr A$.

The arguments in this section are close to those in
\cite[Sect. 4]{ATV}.

The following result is standard.

\begin{lemma}
\label{lem.j.torsion}
If $R$ is a prime noetherian $k$-algebra of finite GK-dimension
and $N \in \mod R$, the following 
are equivalent:
\begin{enumerate}
\item{}
$\GKdim N = \GKdim R$;
\item{}
$j(N)=0$;
\item{}
$N \otimes_R Q \ne 0$, where $Q=\Fract R$.
\end{enumerate}
\end{lemma}

Because $A$ is Gorenstein its dualizing module $\omega_A$ is
invertible, hence isomorphic to $A(\ell)$ for some integer
$\ell$ as a left and as a right module.
Our arguments in this section involve an examination of
the convergent spectral sequence
\begin{equation}
\label{eq.sp.seq}
E_2^{p,q}=\uExt_A^p(\uExt_A^{-q}(M,\omega_A),\omega_A) 
\Rightarrow \HH^{p+q}(M)=
	\begin{cases}
	M & \text{if $p+q=0$,} \\
	0 & \text{if $p+q \ne 0$.}
	\end{cases}
\end{equation}
We will often omit the subscript from $E_2^{pq}$.

\begin{theorem}
\label{thm.3dimAS.E2}
Let $A$ be as above and let $M \in \gr A$. The $E_2$-page of 
the spectral sequence (\ref{eq.sp.seq}) looks like
\begin{equation}
\label{eq.E2.page}
\begin{matrix}
E^{00} & E^{10} & 0 & 0 \\
0 & E^{1,-1} & E^{2,-1} & E^{3,-1} \\
0 & 0 & E^{2,-2} & E^{3,-2} \\
0 & 0 & 0 & E^{3,-3} \\
\end{matrix}
\end{equation}
\end{theorem}
\begin{pf}
Since $\uExt^i_A(-,A)=0$ for $i>3$, the non-zero terms 
on the $E_2$-page of the double-Ext spectral sequence lie in the 
$4 \times 4$ region depicted. Therefore the $E^{20}$ and
$E^{30}$ terms survive to the $E_\infty$-page. But any
non-zero terms on the $E_\infty$ page must lie on the diagonal,
so $E^{20}=E^{30}=0$.
Now $\uExt^3_A(M,\omega_A) \cong H^0_{\fm}(M)^*$ is finite dimensional 
so is Cohen-Macaulay of depth zero, whence
$E^{p3}=0$ for $p<3$. 
This explains the zeroes in the top and bottom rows of
(\ref{eq.E2.page}).

The division ring of fractions, $Q=\Fract A$, 
is flat as a left and as a right $A$-module and
$\gldim Q=0$ so, for $i>0$, 
$$
0=\Ext_Q^i(M \otimes_A Q,A \otimes_A Q) \cong 
Q \otimes_A \uExt^i_A(M,A) .
$$
Applying Lemma \ref{lem.j.torsion} 
to $N=\uExt^i_A(M,A)$, we obtain $E^{0,-1}=E^{0,-2}=E^{0,-3}=0$, giving
the zeroes in the left-most column of (\ref{eq.E2.page}).

It remains to show that $E^{1,-2}=0$.
Set $L=\uExt_A^2(M,A)$; if $L=0$ there is nothing to do, so suppose
that $L \ne 0$. 
Since $\ker(E^{1,-2} \to E^{3,-3})$ survives to the $E_\infty$ page, it
is zero. Since $E^{3,-3}$ 
is finite dimensional so is $E^{1,-2}=
\uExt^1_A(L,A)<\infty$.
If $\uExt^1_A(L,A)=0$ we are finished, so suppose otherwise. 

Consider the $E_2$-page of the spectral sequence for $L$. 
Since $Q \otimes_A L  \cong \uExt_Q^2(M \otimes_AQ,A \otimes_A Q)=0$, 
$\uHom_A(L,A)$ is zero; hence 
the $q=0$ and $q=-1$ rows look like
$$
\begin{matrix}
0 & 0 & 0 & 0 \\
E^{0,-1} & E^{1,-1} & E^{2,-1} & E^{3,-1} \\
\end{matrix}
$$
Since $\uExt_A^1(L,A)$ is non-zero and finite dimensional, 
$\uExt^3_A(\uExt^1_A(L,A),A) \ne 0$.
But the $E^{3,-1}$ term 
survives to the $E_\infty$ page, so must be zero.
From this contradiction we conclude that 
$\uExt^1_A(L,A)=0$, as required.
\end{pf}

\begin{lemma}
\label{lem.even.sum}
If $M \in \gr S$ is a Cohen-Macaulay module, then 
$$
\depth M + \GKdim M \equiv 0 \; (\mod \; 2).
$$
\end{lemma}
\begin{pf}
The lemma is true for any connected graded Gorenstein algebra 
$S$ of finite global dimension, $n$ say, having Hilbert series of
the form $f(t)(1-t)^{-n}$ where $f(t) \in \ZZ[t]$. 
The functional equation \cite[(2.35)]{ATV} relating the Hilbert 
series of a module to that of its dual becomes
$$
H_{M^\vee}(t)=(-1)^dH_M(t^{-1})
$$
when $M$ is a Cohen-Macaulay module of depth $d$ and
$M^\vee=\uExt^{n-d}_S(M,\omega_S)$. If $\GKdim M=r$, then
$H_M(t)=g(t)(1-t)^{-r}$ for some $g(t) \in \ZZ[t,t^{-1}]$ so 
$H_{M^\vee}(t)=(-1)^{d+r}t^rg(t^{-1})(1-t)^{-r}$. However, 
$$
\lim_{t\uparrow 1} \, H_{M}(t)
$$
is positive if $M \ne 0$ and the same applies to
$H_{M^\vee}(t)$, so $d+r$ is even.
\end{pf}

\begin{lemma}
\label{lem.j=1}
Let $A$ be as above and $M \in \gr A$. Suppose  $j(M)=1$ 
and write $M^\vee:=\uExt^1_A(M,A)$. Then $M^\vee$ is Cohen-Macaulay
of depth two and GK-dimension two.
\end{lemma}
\begin{pf}
Since $j(M)=1$, the $E_2$-page of the spectral sequence for $M$ looks
like
$$
\begin{matrix}
0 & 0 & 0 & 0 \\
0 & E^{1,-1} & E^{2,-1} & E^{3,-1} \\
0 & 0 & E^{2,-2} & E^{3,-2} \\
0 & 0 & 0 & E^{3,-3} \\
\end{matrix}.
$$
Both $E^{2,-1}$ and $E^{3,-1}$ survive to the 
$E_\infty$-page, so must be zero. Hence $M^\vee$ is
Cohen-Macaulay of depth two and, by Lemma \ref{lem.even.sum}, $\GKdim
M^\vee=2$.
\end{pf}

\begin{theorem}
\label{thm.3dimAS.CM}
If $A$ is as above, then
\begin{equation}
\label{eq.num.CM}
j(M) +\GKdim M=3
\end{equation}
for all non-zero finitely generated graded $A$-modules $M$.
\end{theorem}
\begin{pf}
By Lemma \ref{lem.j.torsion}, (\ref{eq.num.CM}) holds when
$M=\GKdim 3$ and when $j(M)=0$.
Since finite dimensional modules are precisely the Cohen-Macaulay
modules of depth zero, (\ref{eq.num.CM}) holds when $\GKdim M=0$
too, so it remains to prove
(\ref{eq.num.CM}) for modules of GK-dimensions 1 and 2.

Suppose  $\GKdim M=1$. By Lemma \ref{lem.j.torsion}, $j(M)\ge
1$. As in the proof of Lemma \ref{lem.j=1},
the $E_2$-page of the double-Ext spectral sequence for $M$ 
looks like
$$
\begin{matrix}
0 & 0 & 0 & 0 \\
0 & E^{1,-1} & 0 & 0 \\
0 & 0 & E^{2,-2} & E^{3,-2} \\
0 & 0 & 0 & E^{3,-3} \\
\end{matrix}
$$
Hence the filtration induced on $M$ by the spectral sequence looks like
$M=F^0M=F^1M \supset F^2M \supset \cdots$, so there is a surjective map
$$
M \to F^1M/F^2M = E^{1,-1}_\infty = \ker(E^{1,-1} \to E^{3,-2}).
$$
This gives an exact sequence $M \to E^{1,-1} \to
E^{3,-2}$. Since $\dim_k(E^{3,-2})<\infty$, it follows that
$\GKdim(E^{1,-1}) \le 1$. If $E^{1,-1} \ne 0$, then
$j(\uExt^1_A(M,A))=1$ so, by Lemma  \ref{lem.j=1}, $\GKdim
\uExt^1_A(M,A)^\vee =2$; that is, $\GKdim(E^{1,-1})=2$ , which
contradicts the foregoing. So we must have $E^{1,-1}=0$. Hence
$E_2^{*,-1}=0$, so $\uExt^1_A(M,A)=0$
and $j(M) \ge 2$. However, $j(M) \ne 3$ because $\dim_k M
=\infty$, so $j(M)=2$.

Now suppose that $\GKdim M=2$. By the first paragraph of the
proof $j(M)$ is either 1 or 2. 
Suppose  $j(M)=2$; we seek a contradiction.
Let $\tau M$ be the sum of all finite dimensional graded submodules
of $M$, and consider the exact sequence $0 \to \tau M \to M \to
\Mol \to 0$. It follows easily (cf. \cite[Prop. 2.46]{ATV}) that
$\uExt^i_A(\Mol,A)=0$ for $i \ne 2$ so $\Mol$ is Cohen-Macaulay of
depth one. But $\GKdim \Mol=\GKdim M=2$
which contradicts Lemma \ref{lem.even.sum}. 
\end{pf}

\begin{theorem}
\label{thm.ausl}
The algebra $A$ satisfies the Auslander condition: if $M$ is a finitely
generated $A$-module and $N$ an $A$-submodule of $\uExt^j_A(M,A)$,
then $\uExt^i_A(N,A)=0$ for $i<j$.
\end{theorem}
\begin{pf}
Let $M$ be a non-zero finitely generated graded $A$-module, 
and $N$ a non-zero graded $A$-submodule of 
$\uExt_A^i(M,A)$. By Theorem \ref{thm.3dimAS.E2}, the $E_2$ page of
the spectral sequence for $M$ looks like
(\ref{eq.E2.page}), so $j(\uExt_A^i(M,A)) \ge i$. By Theorem
\ref{thm.3dimAS.CM}, $\GKdim(\uExt^i_A(M,A)) \le 3-i$ so
$\GKdim N \le 3-i$; by Theorem \ref{thm.3dimAS.CM} applied to
$N$, $j(N) \ge i$.
\end{pf}

\begin{corollary}
\label{cor.ausl}
The algebra $S$ satisfies the Auslander condition.
\end{corollary}
\begin{pf}
This follows from \cite[Thm. 3.6]{Lev}.
\end{pf}

Because $S$ satisfies the Auslander condition, the results in
sections 1 and 2 of \cite{LS} apply. In \cite[Sect. 1]{LS} a 
non-zero module $M \in \gr S$ is said to be Cohen-Macaulay if its
projective dimension is the smallest $i$ such that $\Ext_S^i(M,S)$ 
is non-zero. Since $S$ is Gorenstein, $M$ is 
Cohen-Macaulay in the sense of \cite{LS} if and only if it 
is Cohen-Macaulay in the sense of the present paper.

\section{Lines and rulings}
\label{sect.lines}

We continue to assume that $S$ and $A$ are as in the notation just
before Proposition \ref{prop.C(A)}. We continue to use the notation 
$Q:=\Proj A$.

\subsection{Line modules and maximal Cohen-Macaulay modules}

A {\sf graded line module} for $A$ or $S$ is a graded module $L$ that 
is cyclic and has Hilbert series 
$$
H_L(t)=(1-t)^{-2}.
$$
We will write $\cO_L$ for the image of $L$ in either $\Proj A$ or
$\Proj S$. The class of $\cO_L$ in $K_0(\Proj S)$ is $(1-t)^2$.
By way of comparison, if $0 \ne x \in S_1$ and $\cO_H$ denotes the
image of $S/xS$ in $\Proj S$, the class of $\cO_H$ is $1-t$.

\begin{lemma}
\label{lem.mcm}
Let $M \in \MM$.
Then 
\begin{enumerate}
\item{}
$\Hom_{\Gr A}(M(-1),A) \cong k^2$;
\item{}
if $0 \ne f \in \Hom_{\Gr A}(M(-1),A)$ then $f$ is injective, and
\item{}
$\coker f$ is a line module.
\end{enumerate}
\end{lemma}
\begin{pf}
(1)
There is an exact sequence 
$$
0 \to \uHom_A(M,A) \to A^2 \to 
 \uHom_A(\Omega M,A) \to 0
$$
of maximal Cohen-Macaulay {\it left} modules. 
Now $\uHom_A(\Omega M,A)$ is indecomposable because $M$ is, and is
obviously generated by its degree zero component which is
two-dimensional because $\Hom_{\Gr A}(M,A)=0$. Hence
$\uHom_A(\Omega M,A)$ belongs to $\MM'$, the corresponding set of
maximal Cohen-Macaulay {\it left} $A$-modules.
By the left module version of Lemma \ref{lem.two.mcms},
$\uHom_A(M,A)(1)$ is also in $\MM'$, so $\Hom_{\Gr A}(M(-1),A)\cong
\uHom_A(M,A)(1)_0 \cong k^2$.

(2) and (3)
Because $M$ is 3-critical and $A$ is a domain, and hence
3-critical, every non-zero map $M(-1) \to A$ is injective. 
A Hilbert series computation shows that $\coker f$ is a line module.
\end{pf}

\begin{proposition}
\label{prop.lines}
Suppose $S$ is a connected, graded, Gorenstein, Koszul algebra
with Hilbert series $(1-t)^{-4}$. 
Let $z$ be a central regular element of degree two in $S$ and set
$A=S/(z)$. 
Let $L$ be a line module for $A$.  Then
\begin{enumerate}
\item{}
$L$ has a linear resolution as an $S$-module;
\item{}
$L$ has a linear resolution as an $A$-module, namely $\cdots \to
A(-2)^2 \to A(-1)^2 \to A \to L \to 0$;
\item{}
$L$ is Cohen-Macaulay of depth two and 2-critical with respect to
GK-dimension;
\item{}
there is an exact sequence $0 \to M(-1) \to A \to L \to 0$ for 
a unique $M \in \MM$;
\item{}
if $M$ and $L$ are as in part (4), then $FM \cong (FL)(1)$.
\end{enumerate}
\end{proposition}
\begin{pf}
(1)
By \cite[Cor. 2.9]{LS}, the minimal resolution of $L$ 
over $S$ is 
$$
0 \to S(-2) \to S(-1)^2 \to S \to L \to 0.
$$

(2)
If $L$ is {\it any} $A$-module having a linear resolution over $S$,
then $L$ has a linear resolution over $A$: to see this, 
use the fact that $\uExt^1_S(A,k) \cong k(2)$ and use the long 
exact sequence associated to the degenerate spectral sequence 
$$
\uExt^p_A(L,\uExt^q_S(A,k)) \Rightarrow \uExt^{p+q}_S(L,k).
$$
This general fact for commutative rings is proved in \cite{Eis}.

(3)
A line module is Cohen-Macaulay of depth two
by \cite[Prop. 2.8]{LS}, and 2-critical by \cite[Cor. 1.11]{LS}.

(4)
If $M(-1)$ is the kernel of a surjective map $A \to
L$, then it follows from the long exact sequence for local cohomology 
that $M(-1)$ is Cohen-Macaulay of depth three. Every submodule of $A$
is indecomposable because $A$ is a domain. Hence $M(-1)$ is
indecomposable. 
From the linear resolution of $L$, we see that $M(-1)$ is generated by
$M(-1)_1$ and  that $\dim M(-1)_1=2$, so $M \in \MM$.

The uniqueness of $M$ will follow from (5) because if $M'$ is another
element of $\MM$ such that $M' \cong \ker(A \to L)$, then $FM$ is isomorphic to
$FM'$ in $\sD^b(\qgr A^!)$, whence $M$ and $M'$ are isomorphic in $\uMCM(A)$.
Now apply  Lemma \ref{lem.indec.MCM}.

(5)
This follows from Lemma \ref{lem.F.twists} because $M \cong \Omega L(1)$.
\end{pf}

\subsection{Rulings}
For each $M \in \MM$, we define
$$
\hbox{the {\sf ruling corresponding to $M$}}:=
\{L_\phi:=\coker \phi \; | \; \phi \in  \PP(\Hom_{\Gr A}(M(-1),A))\}.
$$
Then
\begin{itemize}
\item{}
each ruling consists of a $\PP^1$ of line modules;
\item{}
every line module belongs to a unique ruling;
\item{}
$Q$ has two rulings if it is smooth, and one otherwise.
\end{itemize}
The first of these facts follows from Lemma \ref{lem.mcm}, 
the second from Proposition \ref{prop.lines}(4),
and the third is a consequence of Theorem \ref{thm.main} and the fact that 
the cardinality of $\MM$ equals the number of isomorphism classes of 
simple $C(A)$-modules.

The results in the rest of this section provide further justification 
for using the word ruling.

\begin{lemma}
\label{lem.lines2}
Let $L_\phi$ and $L_\psi$ $(\phi,\psi \in \PP^1$) be lines in the same ruling.
Then
\begin{enumerate}
\item{}
$L_\phi \cong L_\psi$ if and only if $\phi =\psi$;
\item{}
$\pi L_\phi \cong \pi L_\psi$ if and only if $\phi =\psi$.
\end{enumerate}
\end{lemma}
\begin{pf}
(1)
Let $L_\phi$ and $L_\psi$ be in the ruling corresponding to $M \in \MM$.
Because $M$ is indecomposable, $\Hom_{\Gr A}(M,M)$ is a local ring.
It is finite dimensional and contains no non-zero nilpotents since $M$ is
GK-homogeneous. Hence $\Hom_{\Gr A}(M,M) \cong k$. 
Because $L_\phi$ is cyclic, 
$L_\phi \cong L_\psi$ if and only $\Im \phi =\Im \psi$.
However, the images are the same if and only if $\phi =
\psi \theta$ for some $\theta \in \Hom_{\Gr A}(M(-1),M(-1))$; that is,
if and only if $\phi =\psi$ as elements of $\PP^1=\PP(\Hom_{\Gr A}(M(-1),A)$.

(2)
Because $L_\psi$ is Cohen-Macaulay of depth two,
the exact sequence (\ref{eq.omega.pi})
implies that $\omega\pi L_\psi \cong L_\psi$.
Hence
$$
\Hom_X(\pi L_\phi, \pi L_\psi) \cong \Hom_{\Gr A}(L_\phi,\omega\pi
L_\psi) \cong \Hom_{\Gr A}(L_\phi,L_\psi)
$$
so the result follows from (1).
\end{pf}

The argument in (2) and the observation that each line module belongs to a
unique ruling shows that is $L$ and $L'$ are line modules in different rulings,
 then $\pi L \not\cong \pi L'$.

\begin{proposition}
\label{prop.lines.plane}
Let $L$ and $L'$ be line modules.
\begin{enumerate}
\item{}
If $Q$ is smooth, then
$L$ and $L'$ belong to different rulings if and only if there
is an exact sequence 
\begin{equation}
\label{ses.lines.plane}
0 \to L'(-1) \to A/aA \to L \to 0
\end{equation}
for some $0 \ne a \in A_1$.
\item{}
If $Q$ is not smooth there is always an exact sequence of the form
(\ref{ses.lines.plane}).
\end{enumerate}
\end{proposition}
\begin{pf}
There are exact sequences $0 \to M(-1) \to A \to L \to 0$
and $0 \to M'(-1) \to A \to L' \to 0$ in which $M, M' \in \MM$.

(1)
($\Rightarrow$)
Suppose  $L$ and $L'$ belong to different rulings; then $M
\not\cong M'$ so $M'(-1) \cong \Omega M$ by 
Corollary \ref{cor.M.OmegaM}.
The first term in the exact sequence
$$
%0 \to \Hom(L,L'(-1)) \to 
\Hom_{\Gr A}(A,L'(-1)) \to \Hom_{\Gr A}(M(-1),L'(-1)) \to
\Ext^1_{\Gr A}(L,L'(-1)) \to 0
$$
is zero, and $\Hom_{\Gr A}(M(-1),L'(-1))$ is isomorphic to
\begin{align*}
 \Ext^1_{\Gr A}(M(-1),\Omega M(-1)) & 
\cong \Hom_{\uMCM}(M(-1),\Omega M(-1)[1]) 
\\
& = \Hom_{\uMCM}(M(-1),M(-1)) 
\\
& \cong k.
\end{align*}
Hence $\Ext^1_{\Gr A}(L,L'(-1)) \ne 0$, and
there is a non-split exact sequence 
$$
\begin{CD}
0 @>>> L'(-1) @>>> V @>{\a}>> L @>>> 0
\end{CD}
$$
in $\Gr A$.
Choose $0 \ne \phi \in \Hom_{\Gr A}(A,V)$. The composition
$\a\phi:A \to L$ is surjective because $L$ is cyclic, so $\dim_k
\phi(A_1) \ge 2$. If $\dim_k\phi(A_1)=2$, then $L \cong A/WA$, where
$W=(\ker \phi)_1$, whence the map $V \to L$ splits, contrary to
our assumption. Thus $\dim_k\phi(A_1)=3$. Hence there is some $0
\ne a \in A_1$ and a map $\psi:A/aA \to V$ that is surjective
in degrees zero and one. Let $K=\ker \a\psi$. 
There is a commutative diagram
$$
\begin{CD}
0 @>>> K @>>> A/aA @>>> \im(\a\psi) @>>> 0
\\
@. @VVV @VV{\psi}V @VV{\cong}V
\\
0 @>>> L'(-1) @>>> V @>{\a}>> L @>>> 0
\end{CD}
$$
Since $\dim_k(\im \a\psi)_1< \dim_k (A/aA)_1$, $K_1 \ne 0$, and
hence the map $K \to L'(-1)$ is surjective. It follows that $\psi$
is surjective, and hence injective because $V$ and $A/aA$ have the
same Hilbert series. Hence we get a non-split exact sequence as
claimed.

($\Leftarrow$)
To show that $L'$ is in a different ruling from $L$ it suffices to
show that $\Ext^1_{\Gr A}(L',M(-1))=0$. 

Apply $\Hom_{\Gr A}(-,M(-1))$ to the exact sequence
(\ref{ses.lines.plane}).
A computation shows that $\Ext^1_{\Gr A}(A/aA(1),M(-1))=0$. Since
$\Ext^2_{\Gr A}(A/aA(1),M(-1))=0$, we have $\Ext^1_{\Gr A}(L',M(-1)) \cong
\Ext^2_{\Gr A}(L(1),M(-1))$. 
From the exact sequence $0 \to M \to A(1) \to
L(1) \to 0$, we see that 
$$
\Ext^2_{\Gr A}(L(1),M(-1)) \cong \Ext^1_{\Gr A}(M,M(-1)) \cong 
\Hom_{\uMCM}(M,M(-1)[1])
$$
and this is zero as we see by applying the functor $F$.

(2)
The proof of the implication ($\Rightarrow$) works when $Q$ is
not smooth too because then $M' \cong M \cong \Omega M(1)$.
\end{pf}

\section{Points on quadrics}
\label{sect.pts}

We continue to assume that $S$ and $A$ are as in the notation just
before Proposition \ref{prop.C(A)}.
We also assume $Q$ is smooth.

\subsection{Point modules}
A {\sf graded point module} for $A$ or $S$ is a graded module $P$ that 
is cyclic and has Hilbert series 
$$
H_P(t)=(1-t)^{-1}.
$$
We will write $\cO_P$ for the image of $P$ in either $\Proj A$ or
$\Proj S$. The class of $\cO_P$ in $K_0(\Proj S)$ is $(1-t)^3$.

\begin{lemma}
\label{lem.lines.mcms}
Let $M \in \MM$. 
If $L_\phi$ is in the ruling corresponding to $\Omega M(1)$, 
there is an exact sequence of the form
\begin{equation}
\label{eq.A.M.L}
\begin{CD}
0 @>>> A @>>> M @>>> L_\phi @>>> 0.
\end{CD}
\end{equation}
\end{lemma}
\begin{pf}
By hypothesis, there is an exact sequence
\begin{equation}
\label{eq.M.A.L}
\begin{CD}
0 @>>> \Omega M @>{\phi}>> A @>{\phiol}>> L_\phi @>>> 0,
\end{CD}
\end{equation}
There is also an exact sequence 
\begin{equation}
\begin{CD}
0 @>>> \Omega M @>{\theta}>> A^2 @>>> M @>>> 0.
\end{CD}
\end{equation}
Since $M$ is maximal Cohen-Macaulay $\Ext^1_{\Gr A}(M,A)=0$,
so the natural map
$$
\Hom_{\Gr A}(A^2,A) \to \Hom_{\Gr A}(\Omega M,A) \qquad
\rho \mapsto \rho \circ \theta
$$
is surjective and hence an isomorphism because $\Hom_{\Gr
A}(M,A)=0$. Hence $\phi=\rho \theta$ for a unique $\rho:A^2 \to A$. 

The map $\phiol \rho$ in the diagram
$$
\begin{CD}
0 @>>> \Omega M @>{\theta}>> A^2 @>>> M @>>> 0.
\\
@. @. @VV{\phiol\rho}V
\\
@. @. L_\phi
\end{CD}
$$
is surjective because $\phiol$ and $\rho$ are, and
$\phiol\rho\theta=\phiol\phi=0$, so there is a surjective map
$\psi:M \to L_\phi$.

Now $(\ker \psi)_0 \ne 0$ because $\dim M_0 > \dim (L_\phi)_0$,
so there is a non-zero map $A \to \ker \psi$. Since $M$ is
3-critical so is $\ker \psi$, so the map $A \to \ker\psi$ is
injective. But 
$$
\H_{\ker\psi}(t)=H_{M}(t)-H_{L_\phi}(t)=H_A(t)
$$
so the map $A \to \ker \psi$ must be an isomorphism.
\end{pf}

\begin{proposition}
\label{prop.lines.pts}
Let $M \in \MM$.
If $L_\phi$ is in the ruling corresponding to $\Omega M(1)$,
there is an exact sequence
\begin{equation}
\label{eq.L.L'.P}
\begin{CD}
0 @>>> L(-1) @>>> L_\phi @>>> P @>>> 0
\end{CD}
\end{equation}
in which $P$ is a point module and $L$ is a line module in the same
ruling as $L_\phi$.
\end{proposition}
\begin{pf}
By hypothesis, there is an exact sequence of the form (\ref{eq.M.A.L}). 

\underline{Claim:} 
$\dim \Hom_{\Gr A}(M(-1),L_\phi) \ge 2$.
\underline{Proof:} 
Applying $\Hom_{\Gr A}(M(-1),-)$ to  (\ref{eq.M.A.L}) yields an exact sequence
$$
\begin{CD}
\Hom_{\Gr A}(M(-1),A) @>{\gamma}>> \Hom_{\Gr
A}(M(-1),L_\phi) @>{\delta}>> \Ext^1_{\Gr A}(M(-1),\Omega M).
%\cong \Hom_{\uMCM}(M(-1),\Omega M[1]).
\end{CD}
$$
If $Q$ is smooth, then $\Hom_{\Gr A}(M(-1),\Omega M)=0$
by Corollary \ref{cor.M.OmegaM}, so $\gamma$ is injective, and the 
claim follows from Lemma \ref{lem.mcm}.

Suppose $Q$ is not smooth. 

Then $M(-1) \cong \Omega M$ by Corollary 
\ref{cor.M.OmegaM}; thus $\Hom_{\Gr A}(M(-1),\Omega M) \cong k$ 
by the argument in the proof of Lemma \ref{lem.lines2}, 
and $\Ext^1_{\Gr A}(M(-1),
\Omega M) \ne 0$ by the proof of Step 1 in Theorem \ref{thm.main}.
However,
$\Ext^1_{\Gr A}(M(-1),A)=0$ so the claim holds in this case
too. $\lozenge$

The restriction of each non-zero $\psi \in \Hom_{\Gr A}(M(-1),L_\phi)$
gives a non-zero map $M(-1)_1 \to (L_\phi)_1$ between two 
2-dimensional vector spaces. Now every line
in the projective space $\PP^3=\PP(\Hom_k(k^2,k^2))$
meets the quadric of singular maps, so there is a non-zero $\psi$ such 
that $(\ker\psi)_1 \ne 0$. There is a
non-zero map $A(-1) \to \ker \psi$; this map is injective
because $A(-1)$ and $\ker\psi$ are 3-critical; the cokernel
of the composition $\alpha:A(-1) \to \ker \psi \to M(-1)$ is cyclic
because $\dim M_0 = 1+\dim A(-1)_1$ and $M$ is generated by $M_0$;
the Hilbert series of $\coker\alpha$ is
$$
H_{M(-1)}(t)-H_{A(-1)}(t) = t(1-t)^{-2}
$$
so $\coker \alpha$ is a shifted line module, say $L(-1)$.

Because $\psi\alpha=0$, it follows from the diagram
$$
\begin{CD}
0 @>>> A(-1) @>{\alpha}>> M(-1) @>{\alphaol}>> L(-1) @>>> 0
\\
@. @. @VV{\psi}V
\\
@. @. L_\phi
\end{CD}
$$
that $\psi = \beta\alphaol$ for some $\beta:L(-1) \to L_\phi$.
Because $L(-1)$ and $L_\phi$ are 2-critical $\beta$ is injective, and
$\coker \beta$ is cyclic with Hilbert series 
$(1-t)^{-2}-t(1-t)^{-2}$. Hence $\coker \beta$ is a point module.

It remains only to show that $L$ is in the same ruling as $L_\phi$.
Consider the diagram
$$
\begin{CD}
0@>>> \Omega M(-1) @>>> A(-1)^2 @>{\theta}>> M(-1) @>>> 0.
\\
@. @. @. @VV{\alphaol}V
\\
@. @. @. L(-1)
\end{CD}
$$
Since $L$ is cyclic, the restriction of $\alphaol\theta$
to one of the copies of $A(-1)$ is surjective so, after a Hilbert
series computation, we see that the kernel of $\alphaol\theta$ 
must be isomorphic to $\Omega M(-1)$.
From the exact sequence $0 \to  \Omega M(-1) \to A(-1) \to L(-1) \to 0$
we see that $L$ is in the same ruling as $L_\phi$.
\end{pf}

\begin{lemma}
\label{lem.L.L.P}
Given a line module $L$, there is an exact sequence of the form
\begin{equation}
\label{ses.LLP}
\begin{CD}
0 @>>> L(-1) @>>> L_\phi @>>> P @>>> 0
\end{CD}
\end{equation}
in which $P$ is a point module and $L_\phi$ is a line module in 
the same ruling as $L$.
\end{lemma}
\begin{pf}
This follows from duality and Proposition \ref{prop.lines.pts} for
left modules.

Since $L$ is Cohen-Macaulay of depth two, $\uExt^1_A(L,A)$ is a
{\it left} line module. By Proposition \ref{prop.lines.pts} for
left modules, there is an exact sequence
$$
0 \to L'(-1) \to \uExt^1_A(L,A) \to P' \to 0
$$
where $L'$ and $P'$ are a left line and point module,
respectively. Since $P'$ is Cohen-Macaulay of depth one, applying
$\uHom_A(-,A)$ to this exact sequence gives an exact sequence 
$$
0 \to \uExt^1_A(\uExt^1_A((L,A),A) \to \uExt^1_A(L',A)(1) \to
\uExt^1_A(P',A) \to 0.
$$
Twisting this by $(-1)$ gives the desired exact sequence
(\ref{ses.LLP}).
\end{pf}

\section{The Grothendieck group of a smooth quadric}
\label{sect.K0}

We continue to assume that $S$ and $A$ are as in the notation just
before Proposition \ref{prop.C(A)}.
We also assume that $Q$ is smooth. We will write $C$ for the 
algebra $C(A)$ that is isomorphic to $M_2(k) \oplus M_2(k)$.

We write $K_0(\sA)$ for the Grothendieck group of an abelian category
$\sA$ and define $K_0(Q):=K_0(\coh Q)$. 
We will show that there is an isomorphism 
$K_0(Q) \cong K_0(\PP^1 \times \PP^1)$ of abelian groups 
that is compatible with the Euler forms. However, the discussion
after Theorem \ref{thm.sing.skly.qus} shows that the
effective cones need not coincide under this isomorphism.

\subsection{}
The localization sequence for K-theory gives the exact rows 
in the diagram
\begin{equation}
\label{eq.K0}
\begin{CD}
K_0(\fdim A^!) @>{\a}>> K_0(\gr A^!) @>{\b}>> K_0(\mod C) @>>> 0
\\
@. @V{\phi}V{\cong}V
\\
K_0(\fdim A) @>>> K_0(\gr A) @>>> K_0(Q) @>>> 0.
\end{CD}
\end{equation}
The isomorphism $\phi$ is induced by the Koszul duality
functor $K:\sD^b(\gr A) \to \sD^b(\gr A^!)$ 
and the fact that $K_0(\sA) \cong K_0(\sD^b(\sA))$.
Thus $\phi([N])=[KN]$ is an isomorphism of abelian groups.

We use
the degree shift functors $(\pm 1)$ on the categories $\fdim A^!$,
$\gr A^!$, $\fdim A$, $\gr A$, and $\mod Q$, to make their
Grothendieck groups into $\ZZ[t,t^{-1}]$-modules via
$$
t.[M]=[M(-1)].
$$
Exact functors between these categories that ``commute'' with the
shift functors induce $\ZZ[t,t^{-1}]$-module homomorphisms between
their Grothendieck groups. 

Let $M$ and $M'$ be the two indecomposable maximal Cohen-Macaulay 
modules such that $\MM=\{M(1),M'(1)\}$.

We will use the notation 
$$
a=[A], m=[M], m'=[M'], \ell:=a-m, \ell':=a-m'
$$
for these elements of $K_0(\mod  A)$.
We will use the same notation for the images of $a,m,m',\ell$ 
and $\ell'$ in $K_0(Q)$ and will always take care to indicate which
Grothendieck group we are working in.

The line modules $L$ and $L'$ in the rulings determined by $M$ 
and $M'$ respectively occur in exact sequences of the form
$0 \to M \to A \to L \to 0$ and $0 \to M' \to A \to L' \to 0$,
so all the line modules in a single ruling give the same class 
in $K_0(\mod  A)$, namely $[L]=\ell=a-m$ and $\ell':=[L']=a-m'$
respectively.

\begin{proposition}
\label{prop.K0Q}
The Grothendieck group of $Q$ is free of rank 4 with basis
$\{a,m,m',at\}$ and the action of $\ZZ[t,t^{-1}]$ on it is given by
\begin{align*}
mt & = 2at -m'
\\
m't & = 2at-m,
\\
at^2 & = a(1+4t)-2(m+m').
\end{align*}
As an $R$-module,
\begin{equation}
\label{eq.K0Q}
K_0(Q) \cong {{Ra \oplus R\ell}\over{\bigl(\ell(1-t)^2, \;
a(1-t^2)-2\ell(1-t)\bigr)}}.
\end{equation}
\end{proposition}
\begin{pf}
Write $R=\ZZ[t,t^{-1}]$ and $C=A^![w^{-1}]_0$.  

By D\'evissage, $K_0(\fdim A^!) \cong K_0(\mod  k)$. 
Taking Hilbert series gives a map $K_0(\mod  A^!) \to
\ZZ[[t]][t^{-1}]$; since $ K_0(\mod  k) \cong \ZZ[t,t^{-1}]$
it follows that the map $\a$ in (\ref{eq.K0})
is injective. Since $Q$ is smooth,
$K_0(\mod C) =\ZZ[S_1] \oplus \ZZ[S_2] \cong \ZZ^2$ where $S_1$ and $S_2$
are the two simple $C$-modules; hence the top row of (\ref{eq.K0})
splits, and 
$$
K_0(\mod  A^!) = R[k] \oplus \ZZ[\tilde S_1] \oplus \ZZ [\tilde S_2]
$$
where $\tilde S_1$ and $\tilde S_2$ are the liftings of $S_1$ and $S_2$
via the functor $-\otimes_C A^![w^{-1}]_{\ge 0}$.
Transferring this to $A$ via Koszul duality, we see that
$$
K_0(\mod  A) = Ra \oplus \ZZ m \oplus \ZZ m'.
$$

From the exact sequences in 
Corollary \ref{cor.M.OmegaM.2}, we obtain
relations
\begin{equation}
\label{K0.relns}
2at= mt+m'=m't+m
\end{equation}
in $K_0(\mod  A)$. 
Hence there is a surjective map
\begin{equation}
\label{K0.MCM}
{{Ra \oplus Rm}\over{(2at-m-(2at-mt)t)}} \rightarrow K_0(\mod  A) 
\end{equation}
of $\ZZ[t,t^{-1}]$-modules.
However, it follows from (\ref{K0.relns}) that there is also a
surjective map
$$
K_0(\mod  A) = Ra \oplus \ZZ m \oplus \ZZ m' \to {{Ra \oplus
Rm}\over{(2at-m-(2at-mt)t)}} 
$$
so we conclude that (\ref{K0.MCM}) is an isomorphism.

Since $\ell=a-m$,  we also have
$$
K_0(\mod  A) \cong {{Ra \oplus R\ell}\over{(a(1-t)^2-\ell(1-t^2))}}.
$$
Now $K_0(\fdim A) \cong K_0(\mod  k) \cong \ZZ[t,t^{-1}]$ with basis
$[k] \leftrightarrow 1$, so
$$
K_0(Q) \cong K_0(\mod  A)/([k]),
$$
where $([k])$ denotes the $\ZZ[t,t^{-1}]$-submodule generated by 
$[k]$.

We now compute $[k]$.
From the Hilbert series for $A^!$, we see that
the truncated minimal resolution of $k$ looks like
$$
0 \to N \to A(-2)^7 \to A(-1)^4 \to A \to k \to 0.
$$
It is clear that $N$ is a maximal Cohen-Macaulay module and 
that $N(3)$ has a linear resolution.

Let $F:\mod  A \to \sD^b(\qgr A^!)$ be the functor in Lemma
\ref{lem.F.twists}. Since $N \cong \Omega^3 k$, that lemma shows that
$$
F(N(3)) \cong (Fk)(3) \cong A^!(3).
$$
The equivalence $\qgr A^! \to \mod C$ sends $A^!$ to $C$. 
The degree twist
$(1)$ on $\qgr A^!$ induces an auto-equivalence of $\mod C$,
but every auto-equivalence of $\mod C$ sends ${}_C C$ to ${}_C C$.
Thus, if $G$ is the composition 
$$
\begin{CD}
\mod  A @>F>> \sD^b(\qgr A^!) @>{\sim}>> \sD^b(\mod C),
\end{CD}
$$
then 
$$
G(N(3)) \cong C(3) \cong C.
$$
The functor $G$ sends the two maximal Cohen-Macaulay modules $M(1)$
and $M'(1)$ to the two simple left $C$-modules, so we see that
$$
G(N(3)) \cong G(M(1)^{\oplus 2} \oplus M'(1)^{\oplus 2}).
$$
It now follows from Buchweitz's duality (Theorem \ref{thm.Buch})
that
$$
N \cong M(-2)^{\oplus 2} \oplus M'(-2)^{\oplus 2}.
$$

The truncated resolution of $k$ and equation (\ref{K0.relns}) 
therefore give 
$$
[k]=(1-4t+7t^2)a-(2m+2m')t^2 =(1+4t-t^2)a-2(m+m').
$$
Hence, in $K_0(Q)$, $at^2=a(1+4t)-2(m+m')$.  
It follows that $K_0(Q)$ has basis $\{a,m,m',at\}$, as claimed.
\end{pf}

By Proposition \ref{prop.lines.pts}, $A$ has some graded point 
modules so we define
$$
p:=\ell(1-t) \qquad \hbox{and} \qquad p':=\ell'(1-t)
$$
for the corresponding classes  in $K_0(Q)$.
By (\ref{K0.relns}) that $(m-m')(1-t)=0$, so
$$
p=\ell(1-t)=(a-m)(1-t)=(a-m')(1-t)=\ell'(1-t)=p'.
$$

\begin{proposition}
The sets $\{a,m,m',p\}$ 
and $\{a,\ell,\ell',p\}$
provide $\ZZ$-bases for $K_0(Q)$. 
The $\ZZ[t,t^{-1}]$-action is given by 
$$
a(1-t)=\ell+\ell't=\ell'+\ell t, \quad
\ell(1-t)=\ell'(1-t)=p, \quad p(1-t)=0. 
$$
\end{proposition}
\begin{pf}
Recall that $\ell=a-m$ and $\ell'=a-m'$, so $\{a,\ell,\ell',at\}$
is a basis for $K_0(Q)$.
Furthermore
$$
p=\ell(1-t)=(a-m)(1-t)=a-at-m+(2at-m')=a+at-m-m',
$$
and it follows from this that the two
claimed bases are indeed bases for $K_0(Q)$.
The action of $t$ is already implicit, if not explicit, in the
calculations made in the proof of Proposition \ref{prop.K0Q}.
\end{pf}

The annihilator of $K_0(Q)$ as a $\ZZ[t,t^{-1}]$-module is 
$(1-t)^3$. The submodule of $K_0(Q)$ annihilated by $(1-t)$ is 
$\ZZ p \oplus \ZZ(\ell-\ell')$. 

Taking Hilbert series gives a $\ZZ[t,t^{-1}]$-module
homomorphism $K_0(\mod  A) \to \ZZ[t,t^{-1},(1-t)^{-1}]$,
$[N] \mapsto H_N(t)$. Likewise there is a homomorphism
$q:K_0(\mod  A) \to \ZZ[t,t^{-1}]$ defined by
$$
q[N]=H_N(t)(1-t)^3.
$$
Because $q[k]=(1-t)^3$, there is an induced $\ZZ[t,t^{-1}]$-module
homomorphism 
$$
\qol :K_0(Q) \to \ZZ[t,t^{-1}]/(1-t)^3.
$$
One has
$$
\qol(a)=1+t, \quad \qol(\ell)=\qol(\ell')=1-t, \quad \qol(p)=(1-t)^2.
$$

Suppose  $N \in \mod  A$ has GK-dimension one. Then
$H_N(t)=f(t)(1-t)^{-1}$ for some $f(t) \in \ZZ[t,t^{-1}]$, so
$\qol[\pi N]$ belongs to the ideal of $\ZZ[t,t^{-1}]/(1-t)^3$
generated by $(1-t)^2$. It follows that 
$$
[\pi N] \in \ZZ p \oplus \ZZ(\ell-\ell').
$$

\subsection{The Euler form}
The Euler form on $K_0(Q)$ is denoted by $(-,-)$
and is defined by
$$
([M],[N])=\sum_{i=0}^2 (-1)^i \dim_k \Ext^i_Q(M,N).
$$

\begin{proposition}
\label{prop.euler}
The Euler form on $K_0(Q)$ is given by
$$
(a,a)=(a,\ell)=(a,\ell')=(a,p)=(p,a)=1; \quad
(\ell,a)=(\ell',a)=-1, 
$$
and
$$
(\ell,\ell)=(\ell',\ell')=(\ell,p)=(\ell',p)=(p,\ell)=(p,\ell')=(p,p)=0
$$
and
$$
(\ell,\ell')=(\ell',\ell)=-1.
$$
\end{proposition}
\begin{pf}
Let $P$ be a graded point module occuring in an exact sequence of the form 
$0 \to L_\phi(-1) \to L_\psi \to P \to 0$ where $L_\psi$ and $L_\phi$ 
are line modules in the same ruling.
From the Cohen-Macaulayness of $A,M,M',P$ we see that
$$
(a,a)=1, \quad (a,m)=(a,m')=0, \quad (a,p)=1,
$$
whence
$$
(a,a)=(a,\ell)=(a,\ell')=(a,p)=1.
$$
Serre duality on $Q$ takes the form $\Ext_Q^i(\cF,\cG) \cong
\Ext^{2-i}_Q(\cG,\cF(-2))^*$ for $\cF,\cG \in \mod Q$. Hence
$(x,y)=(y,xt^2)$ for all $x,y \in K_0(Q)$. Also, $(xt,yt)=(x,y)$.

We have $(\ell,a)=(a,\ell t^2)=(a,\ell -2 p)=-1$, and similarly,
$(\ell',a)=-1$. Also, $(p,a)=(a,pt^2)=(a,p)=1$.
In summary,
$$
(\ell,a)=(\ell',a)=-1, \quad (p,a)=1.
$$

Now we show that $(m,m)=1$. The first step is to show that
$\Ext_Q^1(\cM,\cM)=0$. If $0 \to \cM \to \cF \to \cM \to 0$ is
exact, then applying $\omega$ gives an exact sequence $0 \to M \to
F \to M \to 0$ because $R^1\omega M=H^2_{\fm}(M)=0$ and $\omega\pi
M \cong M$. But $\Ext_{\Gr A}^1(M,M) \cong \Hom_{\uMCM}(M,M[1])=0$,
where the last equality follows by applying the functor $F$, so the
sequence in $\Gr A$ splits; but the original sequence in $\Qcoh
Q$ is obtained by applying $\omega$ to this split sequence, so it
splits too. Hence $\Ext_Q^1(\cM,\cM)=0$. Now, $\Ext^2_Q(\cM,\cM)
\cong \Hom_Q(\cM,\cM(-2))^* \cong \Hom_{\Gr A}(M,M(-2))^* =0$
because $M(-2)_1=0$. Finally, $\Hom_Q(\cM,\cM) \cong \Hom_{\Gr
A}(M,M)=k$, so $(m,m)=1$.

Using this gives
$$
(\ell,\ell)=(a-m,a-m)=(a,a-m)-(m,a) +(m,m) = 1 -(a-\ell,a)+1=0,
$$
and similarly, $(\ell',\ell')=0$.
Using Serre duality, we obtain
$$
0=(\ell,\ell)=(\ell,\ell t^2)=(\ell,\ell-2p) =-2(\ell,p)
$$
which gives $(\ell,p)=0$.
Therefore
$$
(\ell,\ell')=(\ell, \ell' + \ell -p)=(\ell,\ell'+\ell t)
=(\ell,a(1-t)) = -1 -(\ell,at) 
$$
and hence
$$
(\ell,\ell')=-1-(at,\ell t^2)=-1-(a,\ell t) =-1 -(a,\ell -p)=-1.
$$
Similarly, $(\ell',\ell)=-1$.

Finally, 
$$
(p,p)=(\ell(1-t),p)=-(\ell t,p) = -(\ell,pt^{-1})= -(\ell,p)=0,
$$
and
$$
(p,\ell)=(\ell,pt^2)=(\ell,p)=0.
$$
This completes the proof.
\end{pf}

Proposition \ref{prop.euler}
is exactly as in the commutative case---of course, our proof
applies to that case too.

\subsection{The intersection pairing on $Q$}
Recall that $\ell$ and $\ell'$ are the classes in $K_0(Q)$ of $\cO_L$ and $cO_{L'}$ where $L$ and $L'$ 
are line modules belonging to different rulings.

The interpretation of the equality  $(\ell,\ell')=-1$ in Proposition \ref{prop.euler}
is that a line in one ruling meets a line in the other ruling with
multiplicity one. In the commutative case this means the two
lines span a hyperplane. Proposition \ref{prop.lines.plane}(1) is 
the appropriate analogue of this.
It is therefore sensible to introduce the notation
$$
h:=[\cO_Q]-[\cO_Q(-1)]=\ell+\ell' t=\ell'+\ell t.
$$

We now define an intersection pairing on $K_0(Q)$ by 
$$
b.c:=-(b,c).
$$

The next calculation shows that everything behaves as it does for points and lines on a smooth quadric surface in $\PP^3$, i.e., as for $\PP^1 \times \PP^1$. 

\begin{proposition}
The intersection pairing has the following properties:
\begin{align*}
\ell.\ell = \ell.p=h.p=p.h=p.\ell=\ell'.\ell' &=0;
\\
\ell.\ell'=\ell.h=\ell'.h=h.\ell'=h.\ell &=1;
\\
\ell.\ell'=\ell'.\ell &=1.
\end{align*}
\end{proposition}
\begin{pf}
The calculations are as follows:
\begin{align*}
(h,\ell) & =(\ell+\ell' t,\ell)=(\ell't,\ell)=(\ell'-p,\ell)=-1;
\\
(h,p)& =(\ell +\ell't,p)=(\ell,p)+(\ell',pt^{-1})=0;
\\
(\ell,h)& =(h,\ell t^2) = (h,\ell-2p)=-1;
\end{align*}
and $
(p,h)=(h,pt^2)=(h,p)=0.$
\end{pf}

One other computation of interest is $(m,m')=0$.

\begin{proposition}
\label{prop.lines.ext}
Suppose $Q$ is smooth.
Let $L$ and $L'$ be non-isomorphic line modules for $A$, and 
$\cO_L$ and $\cO_{L'}$ their images in $\Proj Q$. The following are
equivalent:
\begin{enumerate} 
\item{}
$\Ext_Q^1(\cO_L,\cO_{L'})=0$;
\item{}
$\Ext^i_Q(\cO_L,\cO_{L'})=0$ for all $i$;
\item{}
$L$ and $L'$ belong to the same ruling.
\end{enumerate}
\end{proposition}
\begin{pf}
Because $L'$ is Cohen-Macaulay of depth 2, $\omega\pi L' \cong L'$.
Hence 
$$
\Hom_Q(\cO_L,\cO_{L'}) \cong \Hom_{\Gr A}(L,\omega \pi L') \cong
\Hom_{\Gr A}(L,L')=0.
$$
By Serre duality 
$\Ext^2_Q(\cO_L,\cO_{L'}) \cong \Hom_Q(\cO_{L'},\cO_L(-2))$;
because $L$ is Cohen-Macaulay of depth 2, this is isomorphic to
$\Hom_{\Gr A}(L',L(-2))$ which is zero because $L(-2)_0=0$.
Hence (1) $\Leftrightarrow$ (2).

By Proposition \ref{prop.euler}, $L$ and $L'$
belong to the same ruling if and only if $([\cO_L],[\cO_{L'}])=0$.
This, together with the observations in the previous paragraph,
shows that (3) is equivalent to (1) and (2).
\end{pf}

If $Q$ is not smooth then $\Ext^1_Q(\cO_L,\cO_{L'}) \ne 0$ and, if
$L \not \cong L'$, then $\Ext^2_Q(\cO_L,\cO_{L'})$ and
$\Hom_Q(\cO_L,\cO_{L'})$ are both zero.

\section{The Sklyanin quadrics}
\label{sect.skly}

Throughout this section $S$ denotes a four-dimensional Sklyanin algebra and  
$$
\PP^3_{\sf Skly} = \Proj S.
$$
We recall some results from 
\cite{LS}, \cite{Sm1}, \cite{Sm4}, \cite{SS}, \cite{TV}, and \cite{vdB2}.

\subsection{}
The data used to define $S$ is a triple $(E,\cL,\tau)$ consisting 
of an elliptic curve $E$, a degree four line bundle $\cL$ on it, 
a translation automorphism $\tau$ of $E$, and
$S$ is a quotient of the tensor algebra on $H^0(E,\cL)$ having Hilbert series $(1-t)^{-4}$.  Like the polynomial ring,
$S$ is Gorenstein, and its dualizing module is $\omega_S \cong 
S(-4)$ as a one-sided $S$-module.
Furthermore, $S$ is a noetherian domain and a Koszul algebra. 
Thus $\PP^3_{\sf Skly} $ is a quantum $\PP^3$ in the sense of section
\ref{sect.qu.Pn}.

Because $S_1=H^0(E,\cL)$ we can, and will, consider 
$E$ as a fixed quartic curve in $\PP(S_1^*)$.
We fix an origin $0$ for $E$ such that four points of $E$ are coplanar
if and only if their sum is zero. 
We therefore identify $\tau$ with a point
on $E$ so the translation automorphism becomes $p \mapsto p+\tau$.

\subsection{Pencils of quadrics in $\PP^3$ and $\PP^3_{\sf Skly}$}
A generic pencil of quadrics in $\PP^3$ has exactly four singular members. Its base locus is
a quartic elliptic curve. The smooth quadrics have two rulings on them, and the singular 
ones have only one ruling. The lines on the quadrics are the secant lines to the base locus.

The pencil of commutative quadrics in $\PP(S_1^*)$ containing $E$
may be labelled as $Y_z$, $z \in E/{\pm} \cong \PP^1$, in such a 
way that $Y_z$ is the union of the secant lines $\overline{pq}$ 
such that $p+q=z$. 
It follows that $Y_z=Y_{-z}$ and the four singular quadrics are
$Y_\omega$, $\omega \in E_2$, the 2-torsion subgroup of $E$.
When $z \notin E_2$, the two rulings on $Y_z$ are given by
$\{\overline{pq} \; | \; p+q=z\}$ and
$\{\overline{pq} \; | \; p+q=-z\}$.

As we now explain, the Sklyanin quadrics behave in a similar way.

The center of $S$ contains two linearly independent
homogeneous elements, $\Omega_1$ 
and $\Omega_2$, of degree two. These give rise to a pencil of quotients
$A=S/(\Omega)$, $\Omega$ a non-zero linear combination of 
$\Omega_1$ and $\Omega_2$, and hence a pencil of 
non-commutative quadric hypersurfaces $\Proj A \subset 
\PP^3_{\sf Skly}$. Each $A$ is a Gorenstein domain  with dualizing
module $\omega_A \cong A(-2)$ as a one-sided $A$-module.

Since $S/(\Omega_1,\Omega_2)$ is a twisted homogeneous 
coordinate ring of 
$E$, $\Proj S/(\Omega_1,\Omega_2)$ presents $E$ as a closed subspace
of $\PP^3_{\sf Skly}$. It is the base locus of the pencil of non-commutative
quadrics.

\subsection{}
The following rule sets up a bijection between the line modules for $S$ 
and the secant lines to $E$ in $\PP(S_1^*)$: if $p,q \in E$, and 
$W \subset S_1$ is the subspace of linear forms vanishing on the 
 $\overline{pq}$, then $S/SW$ is a 
line module that we denote by $L(\overline{pq})$  \cite{LS}.

If $z \in E$, there is a non-zero linear combination $\Omega(z)$ 
of $\Omega_1$ and $\Omega_2$ such that
$$
\Omega(z).L(\overline{pq})=
0 \Longleftrightarrow p+q=z \hbox{ or } p+q=-z-2\tau
$$
(see \cite[Sect. 6]{LS}).
We label the non-commutative quadrics in $\PP^3_{\sf Skly}$ as
$$
Q_z:=\Proj S/(\Omega(z)), \qquad z \in E.
$$
Thus $Q_z=Q_{-z-2\tau}$.

\subsection{Families of lines}
If $z \notin E_2+\tau$, we say there are {\sf two families} of line
modules for $A$ giving ``lines'' on $Q_z$, namely  
$\{L(\overline{pq}) \; | \; p+q=z\}$ and
$\{L(\overline{pq}) \; | \; p+q=-z-2\tau\}$.

The degree two divisors $(p)+(q)$ such that $p+q=z$ are 
parametrized by the points in the fiber over $z$ of the addition map 
$S^2E \to E$. These fibers are isomorphic to
$\PP^1$, which is why we say these lines form a family. 

The next result shows that these ``families'' coincide with the
``rulings'' defined in section \ref{sect.lines}.

\begin{proposition}
\label{lines.family}
Let $L$ and $L'$ be line modules for $A$. Then $L$ and $L'$ belong
to the same ruling if and only if they belong to the same family. 
\end{proposition}
\begin{pf}
Let $Q_z=\Proj A$.
Suppose  $L=L(\overline{pq})$ and $L'=L(\overline{p'q'})$ where 
$p+q,p'+q' \in \{z,-z-2\tau\}$.

($\Leftarrow$)
Suppose  $p+q=p'+q'$. There are points $r,s \in E$ such 
that $p$, $q$, $r$, and $s$ span a secant plane, say that given by 
$a=0$ for $0 \ne a \in A_1$, and
$p'$, $q'$, $r$, and $s$ also span a secant plane, say that given by 
$b=0$ for $0 \ne b \in A_1$. 

Set $L''=L(\overline{r-\tau,s-\tau})$.
By the argument in the proof of 
\cite[Lemma 4.5]{SmStan}, there are exact sequences
$$
0 \to L''(-1) \to A/Aa \to L \to 0
$$
and
$$
0 \to L''(-1) \to A/Ab \to L' \to 0.
$$
By Proposition \ref{prop.lines.plane}, $L$ and $L''$ belong to
different rulings, and so do $L'$ and $L''$; hence $L$ and $L'$
belong to the same ruling.

($\Rightarrow$)
Suppose  $p+q \ne p'+q'$. 
In this case $p+q+(p'+\tau)+(q'+\tau)=0$, so $p,q,p'+\tau,q'+\tau$
span a secant plane. By \cite[Lemma 4.5]{SmStan}, there is an
exact sequence of the form $0 \to  L'(-1) \to A/xA \to L \to 0$,
so $L$ and $L'$ belong to the same ruling by  Proposition
\ref{prop.lines.plane}.
\end{pf}

\begin{theorem}
\label{thm.sing.skly.qus}
The Sklyanin quadric $Q_z$ is smooth if and only if $z+\tau \notin
E_2$.
The four singular quadrics are $Q_{\omega -\tau}$, $\omega \in E_2$.
\end{theorem}
\begin{pf}
If $z+\tau \notin E_2$, then $Q_z=Q_{-z-2\tau}$ has two families of 
line modules, namely $L(\overline{pq})$ such that $p+q=z$ and
$p+q=-z-2\tau$, whereas if $z + \tau \in E_2$, 
there is only one family of line modules for $Q_z$, 
namely $L(\overline{pq})$ such that $p+q=z=-z-2\tau$.
Now by Theorem \ref{lines.family}, there are two rulings on $Q_z$
if and only if $z+\tau \notin E_2$, so the result follows from
Theorem \ref{thm.main}.
\end{pf}

\subsection{Singular quadrics in a pencil}
There is one significant way in which the pencil of Sklyanin quadrics differs from a generic
pencil of quadrics in $\PP^3$.

The singular locus of a singular quadric $Q$ belonging to  a generic 
pencil in $\PP^3$ is a point, and that point lies on all the lines 
on $Q$. However, the results in \cite{SmStan} 
(see also \cite[Sect. 10]{Sm1}) show there is no analogous result
for the Sklyanin quadrics. For simplicity, we will explain this only 
when $\tau$ has infinite order. 

When $\tau$ has infinite order the closed points 
in $\PP^3_{\sf Skly}$ consist of those on $E$
and a discrete family that may be labelled as 
$$
\{p_{\omega+i\tau} \; | \; \omega \in E_2, i \in \NN \}
$$
in such a way that
\begin{enumerate}
  \item[(a)]
  $p_{\omega+i\tau}$ lies on the 
non-commutative secant line $\overline{pq}$ if and only if 
$p+q=\omega+i\tau$, and
  \item[(b)]
    if $\cF_{\omega+i\tau} \in \Qcoh Q$ is the simple
module corresponding to $p_{\omega+i\tau}$, then
$\dim_k H^0(\PP^3_{\sf Skly},\cF_{{\omega+i\tau}})=i+1$. (The letter $\cF$ stands for {\it fat point}.)
\end{enumerate}
Thus, all the lines in one of the two rulings on the smooth quadrics 
$Q_{\omega+i\tau}=Q_{\omega-(i+2)\tau}$, 
$i \in \NN$, pass through a common point.
The lines on a singular quadric $Q_{\omega-\tau}$ do 
not pass through a common point.

Let $i \in \NN$ and $\omega \in E_2$.
By \cite[Sect. 4]{SmStan}, if 
$p+q=\omega+i\tau$ there is an exact sequence
\begin{equation}
\label{ses.fat.pts}
0 \to \cO_{\overline{p-(i+1)\tau,q-(i+1)\tau}}(-i) \to
\cO_{\overline{pq}} \to \cF_{{\omega+i\tau}} \to 0
\end{equation}
of $Q_{\omega+i\tau}$-modules;
because $(p-(i+1)\tau)+(q-(i+1)\tau)\ne p+q$,
the two lines in (\ref{ses.fat.pts}) belong to different rulings;
it also follows from (\ref{ses.fat.pts}) that the class of
$\cF_{{\omega+i\tau}}$ in $K_0(Q_{\omega+i\tau})$ is
$$
[\cF_{{\omega+i\tau}}]=\ell-\ell't^{i+1}=\ell-\ell'+(i+1)p.
$$
This shows that the positive cone of $K_0(Q_{\omega+i\tau})$
is not the same as that of $K_0(\PP^1 \times \PP^1)$.
A computation in $K_0(Q_{\omega+i\tau})$ using Proposition
\ref{prop.euler} gives
$$
(\cF_{{\omega+i\tau}},\cF_{{\omega+i\tau}})=2,
$$
so $p_{\omega+i\tau}$ behaves like a curve with self-intersection
$-2$.

\subsection{}
Similar behavior is exhibited by the primitive quotient
rings of the enveloping algebra of $\fsl(2,\CC)$ (cf. \cite{LB},
\cite{JTS82} and \cite{vdB2}). More precisely, the homogenized
enveloping algebra of $\fsl(2,\CC)$ is the coordinate ring of a
quantum $\PP^3$ that contains a pencil of non-commutative quadrics
and those non-commutative quadrics behave like the Sklyanin
quadrics. In particular, the finite-dimensional irreducible
representations of $\fsl(2,\CC)$ provide points on certain of these
quadrics that also behave like $-2$-curves---they have self-intersection $-2$.

 \subsection{}
The quadrics in a generic pencil in $\PP^3$ can be viewed as 
the fibers of a family $X \to \PP^1$. 
The total space $X \subset \PP^3 \times \PP^1$ is smooth. It 
seems likely that the analogous non-commutative 3-fold $X_{nc}
\subset \PP^3_{\sf Skly} \times \PP^1$ is also smooth, but we do not know how to
tackle this problem.

\subsection{}
Our methods  apply to the pencil of non-commutative
quadrics in the non-commutative $\PP^3$ associated to 
the enveloping algebra of $\fsl(2,\CC)$.
This pencil of non-commutative quadrics is analogous to
the commutative pencil of quadrics generated by a double 
plane $w^2=0$ and $x^2+y^2+z^2=0$. 
The non-commutative pencil contains a ``double plane'' and one more
singular non-commutative quadric that corresponds to the 
unique primitive quotient of $U(\fsl(2,\CC))$ having 
infinite global dimension.
That particular quotient of $U(\fsl(2,\CC))$ is a simple ring so
has no finite dimensional simple module; this is
analogous to the fact that the singular Sklyanin quadrics are
not the ones having a point that causes infinite global dimension.
The homological properties of the quotients of $U(\fsl_2)$ are
described in \cite{JTS82}.

\subsection{}
Let $Q$ be a smooth non-commutative quadric surface.
It would be interesting to show that there is a map $Q \to \PP^1$
in the sense of \cite[Defn. 2.3]{SmSub}, to define and study 
the fibers of such a map, and to show that $Q$ is
the disjoint union of these fibers in a suitable sense.
It would also be interesting to examine quadric hypersurfaces in
non-commutative analogues of $\PP^n$ for $n>3$.

\end{document}